\documentclass{amsart}

\usepackage{amsmath,amssymb,amsfonts}
\usepackage{amsthm}
\usepackage{bbm,hyperref}
\usepackage[dvipsnames]{xcolor}
\usepackage{graphicx,tikz}

\def\ZZ{{\mathbb Z}}

\def\vv{{\mathbf v}}

\theoremstyle{plain}
\newtheorem{theorem}{Theorem}[section]

\newtheorem{lemma}[theorem]{Lemma}
\newtheorem{Assumption}[theorem]{Assumption}

\theoremstyle{remark}
\newtheorem{remark}[theorem]{Remark}

\title[]{Asymptotic delta-like behavior of the weighted taper in Thomson's multitaper F-test}

\author{Hau-Tieng Wu}

\address{Courant Institute of Mathematical Sciences, New York University, New York, 10012, USA}

\begin{document}
\maketitle

\begin{abstract}
We study the weighted taper arising in Thomson's multitaper F test, which is constructed from pre-plunge discrete prolate spheroidal wave functions (DPSWF). Numerical evidence suggests that, in the asymptotic regime, the Fourier transform of this taper converges to a delta measure supported at the origin. We provide a rigorous theoretical justification for this behavior. 
\end{abstract}

\section{Introduction}
 
Discrete prolate spheroidal sequences (DPSS) and discrete prolate spheroidal wave functions (DPSWF) have been studied extensively since Slepian's seminal 1978 work \cite{slepian1978pswfV}. See, for example \cite{delsarte1985discrete,karnik2020improvedchardpsseigenval,said2023non}. They are closely related to continuous prolate spheroidal wave functions (PSWFs), which have been extensively studied since the seminal works of Slepian and Pollak in 1961 \cite{slepian1961prolate}. See, for example \cite{slepian1965some,liirosenblatt2008mtvariance,osipov2013certain,osipov2013prolate,bonami2014uniform,bonami2016uniform}, which is far from exhaustive. DPSS and DPSWF have found widespread applications in signal processing and spectral time series analysis \cite{mtthomson,mtwelchcomparison,waldenspectral,mtreview,walden2000multivarmt,luismtbias,karnik2019fast}, particularly in multitaper spectral estimation and oscillation detection using Thomson's multitaper $F$-test. In this paper, we focus on studying the weighted taper that arises in Thomson’s multitaper $F$-statistic.
 
Let $N\in \mathbb{N}$ denote the sample size and $W\in (0,1/2)$ the bandwidth. Let $U_k:[-1/2,1/2]\to \mathbb{C}$ denote the $k$-th DPSWF, and let $\vv_k\in \mathbb{R}^N$ denote the $k$-th DPSS, or taper. 
We assume, without loss of generality, that $U_k(0)$ are real after a conventional phase adjustment \cite{slepian1978pswfV}.
The precise relationship between $U_k$ and the Fourier transform of $\vv_k$, including the proper phase factor, will be carefully reviewed below.
Given a time series $X_0,\ldots, X_{N-1}$, 
define the {\em tapered Fourier transform} associated with the $k$-th taper $\vv_k$ by
\[
\hat X_{\vv_k}(f) := \sum_{j=0}^{N-1} \vv_{k,j}X_je^{-i2\pi fj}
\]
where $f\in [-1/2,1/2]$.
Thomson's multitaper $F$-statistic is defined as
\[
T(f):=\frac{(K-1)|\tilde{\mu}(f)|^2\sum_{k=0}^{K-1}U_k(0)^2}{\sum_{k=0}^{K-1}|\hat X_{\vv_k}(f)-U_k(0)\tilde{\mu}(f)|^2}\,,
\]
where $f\in [-1/2,1/2]$, $K< 2NW$ is the number of tapers used, and
\[
\tilde{\mu}(f):=\frac{\sum_{k=0}^{K-1}U_k(0)\hat X_{\vv_k}(f)}{\sum_{k=0}^{K-1}U_k(0)^2}\,.
\]
Thomson's multitaper $F$-statistic is used to detect deterministic sinusoidal components in a given time series. Its underlying rationale has two key elements. First, $\tilde{\mu}(f)$ captures the contribution of an oscillatory component at frequency $f$. It can be appreciable when $f$ coincides with the oscillation frequency and is approximately zero otherwise. Second, because the tapers are orthonormal, the corresponding tapered Fourier coefficients, $\hat X_{\vv_k}$ and $\hat X_{\vv_{k'}}$ are uncorrelated when $k\neq k'$, under the usual assumptions on the stochastic component of the time series. Averaging across tapers can therefore reduce stochastic variability. This motivates interpreting $|\hat Z_{\vv_k}(f)-U_k(0)\tilde{\mu}(f)|^2$ as a tapered estimate of the noise power spectrum after removing the putative oscillatory component.

Define a new taper $w\in \mathbb{R}^N$ based on the weighted sum of DPSS tapers: 
\begin{align}\label{def: multitaper equivalent window expression}
w_j:=\frac{\sum_{k=0}^{K} U_{k}(0)\vv_{k,j}}{\sqrt{\sum_{k=0}^{K}U_{k}(0)^2}}\,,
\end{align}
where $K<2NW$.
We call $w$ the {\em weighted taper}, which has unit $L^2$ norm $\sum_{j=0}^{N-1}w_j^2=1$ since $\langle\vv_k,\vv_l\rangle=\delta_{k,l}$. Note that the numerator of Thomson's multitaper $F$-test \cite{mtthomson,waldenspectral} can be rewritten via
\[
\Big(\sum_{k=0}^KU_{k}(0)^2\Big)^{1/2}\hat{\mu}(f)=\frac{\sum_{k=0}^KU_{k}(0)\hat{X}_{\vv_{k}}(f)}{\sqrt{\sum_{k=0}^KU_{k}(0)^2}}=\sum_{j=0}^{N-1} w_jX_je^{-i2\pi jf}\,.
\]
 
\begin{figure}[htb!]
\includegraphics[width=\textwidth]{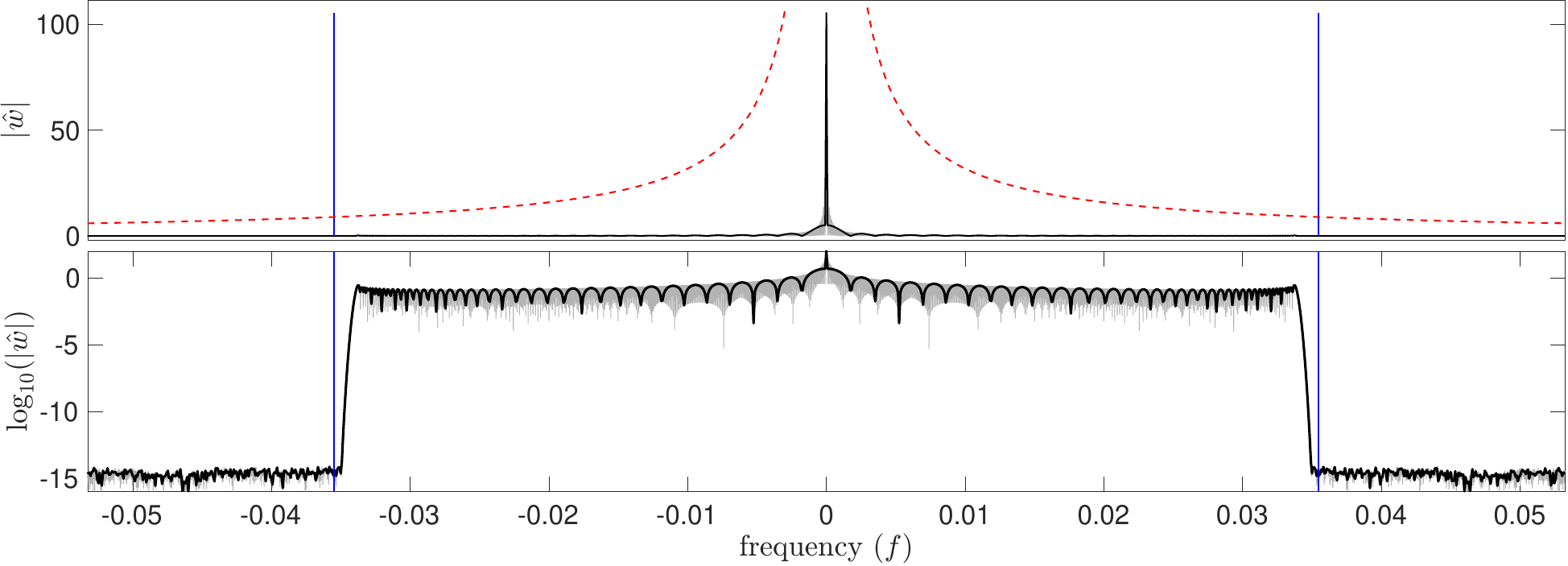}
\includegraphics[width=\textwidth]{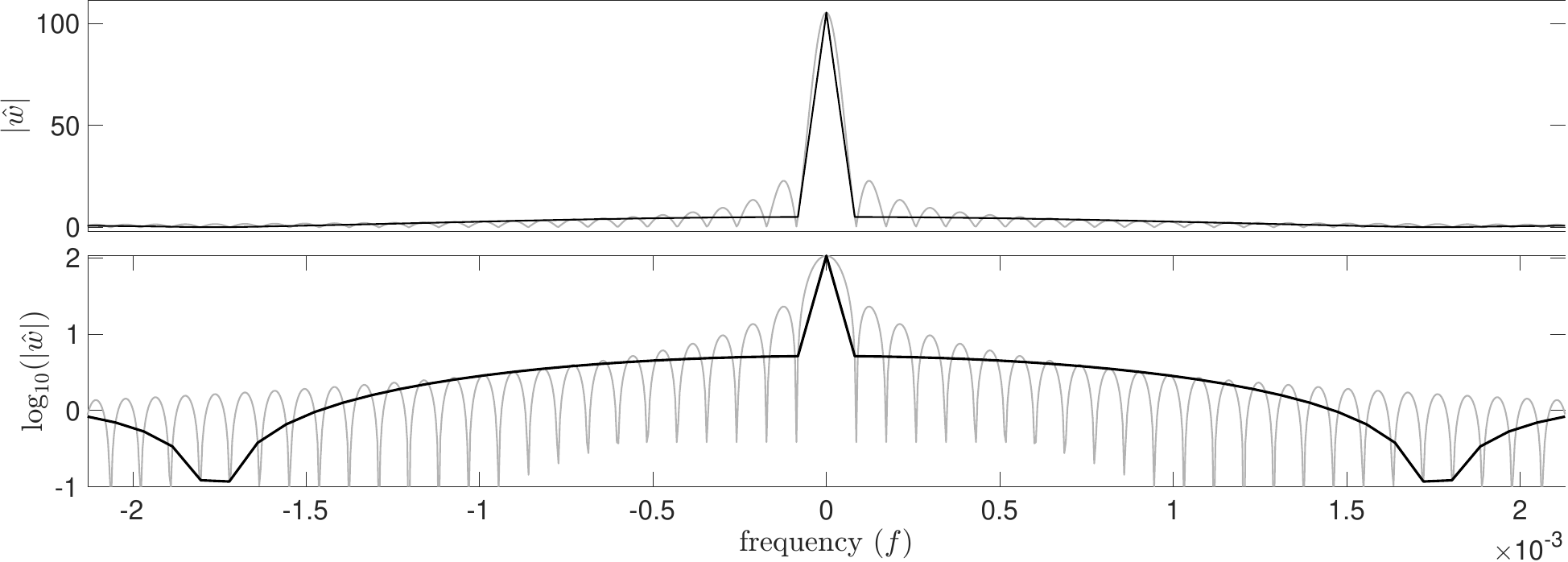}
\caption{First row: $|\hat{w}(i/N)|$, where $i=-(N-1)/2,\ldots,(N-1)/2$, when $N=12,199$, $NW=433$,  $W=0.0355$, and $K=2NW-\lfloor \frac{12}{\pi^2}\log(100NW+25)\log(N)\rfloor$, is shown in black. The red dash is the decay of $1/\sin(\pi\xi)$, $|\hat{w}(f)|$ is shown in gray, and the vertical blue lines are $\pm W$. At frequency $0$, the value is $105.32\approx \sqrt{N}=110.45$. The x-axis is restricted on $[-1.5W,1.5W]$. 
Second row: log scale of the first row; that is, $\log_{10}(|\hat{w}(i/N)|)$ is shown in black and $\log_{10}(|\hat{w}(f)|)$ is shown in gray.  
Third and fourth rows: zoomed in of the first and second rows to $[-0.06W,0.06W]$ for details.}
\label{fig:Vplot}
\end{figure}

Numerically, we found that the discrete time Fourier transform of $w$,
\[
\hat{w}(f)=\sum_{j=0}^{N-1}w_je^{-i2\pi jf}=\frac{\sum_{k=0}^KU_{k}(0)U_{k}(f)}{\sqrt{\sum_{k=0}^KU_{k}(0)^2}},
\] 
where $f\in [-1/2,1/2]$, behaves like a delta-measure supported at $0$. See Figure \ref{fig:Vplot} for an illustration. See \url{https://github.com/hautiengwu2/WeightedTaperThomonsFtest.git} for the Matlab code to reproduce the figure. On the log scale, there is a clear sharp jump at $W$. This numerical finding suggests why in practice Thomson's multitaper $F$-statistic can well detect oscillation with high resolution. However, to our knowledge, this behavior has not been well studied in the literature.

The main contribution of this paper is proving the following theorem to justify the numerical findings. 

\begin{theorem}\label{theorem: bound Vf energy leakage uniform version}
Assume $W=W(N)\to 0$ and $WN\to \infty$ as $N\to \infty$ 
 and $K=\lfloor 2NW\rfloor-G$, where $G:=\lfloor \frac{4c}{\pi^2}\log(100NW+25)\log(N)\rfloor$ with $c>3$. Denote $f_*:=\Big(1- \frac{G}{(NW)^{1/4}}\Big)W$. Take a small constant $\tau>0$. When $N$ is sufficiently large, we have
\begin{itemize}
\item 
$\displaystyle
|\hat{w}(0)|=\sqrt{N}+O\Big(\frac{\log(NW)\log(N)}{N^{1/2}W}\Big)$;
\medskip
\item for $\gamma>1$,
$\displaystyle
\Big|\hat{w}\Big(\frac{1}{\gamma N}\Big)\Big|=\frac{\gamma}{\pi}\sin\Big(\frac{\pi}{\gamma}\Big) \sqrt{N}+O\Big(\frac{\log(NW)\log(N)}{N^{1/2}W}\Big)$;
\medskip
\item   
$\displaystyle
|\hat{w}(k/N)|=O\Big(\frac{\log(NW)\log(N)}{N^{1/2}W}\Big)$
when $k=1,2,\ldots, \lfloor \frac{NW}{\log(NW)}\rfloor$ and
$\displaystyle
|\hat{w}(f)|\leq \frac{\sqrt{N}}{3 k}+O\Big(\frac{\log(NW)\log(N)}{N^{1/2}W}\Big)$
when $f\in \big(\frac{k}{N},\frac{k+1}{N}\big)$;
\medskip
\item 
$\displaystyle
|\hat{w}(f)|= O\Big(\frac{\log(NW)\log(N)}{N^{3/8}W^{7/8}}\Big)$ when $f\in[ \frac{W}{\log(NW)},f_*]$;
\medskip
\item 
$\displaystyle
|\hat{w}(f)|=O\Big(\frac{\log(NW)\log(N)}{N^{3/8}W^{9/8}}\Big)$ when $f\in(f_*,W]$;
\medskip
\item 
$\displaystyle
|\hat{w}(f)|=O(N^{-3/2}W^{-1/2})$ when $f\in (W,(1+\tau)W]$ and $\displaystyle
|\hat{w}(f)|=O(N^{-7/2}W^{-1/2})$ when $f\in \big((1+\tau)W,\frac{1}{2}\big]$.

 \end{itemize}
\end{theorem}

Note that when $W=N^{-\alpha}$ with $\alpha<\frac{1}{3}$, $N^{3/8}W^{9/8}\to \infty$ when $N\to \infty$; that is, despite $\hat{w}(0)$ is of size $\sqrt{N}$ and some sidelobes, when $f\geq \frac{W}{\log(NW)}$, $\hat{w}(f)\to 0$ as $N\to \infty$.  Or if we impose $\alpha<\frac{7}{9}$, we have $\frac{|\hat{w}(f)|}{|\hat{w}(0)|}\to 0$ when $f\geq \frac{W}{\log(NW)}$ as $N\to \infty$. These conditions are sufficiently flexible for practice.
The proof is postponed to Section \ref{section proof of theorem theorem: bound Vf energy leakage uniform version} after reviewing existing results in Section \ref{section Mathematical Background} and establishing some technical lemmas in Section \ref{section technical lemmas}. The basic idea is carefully studying the Sturm-Liouville equation associated with DPSWF, which has one turning point and one degenerate point. The Liouville-Green transform is applied to establish the control of $U_k$ at $0$ and $1/2$, the envelope monotonicity is applied to control $U_k'(1/2)$, and the energy method is applied to control DPSWF over $(W,1/2)$ or $(0,W)$, depending on $k$. Finally, combine recent results reported in \cite{karnik2020improvedchardpsseigenval,said2023non} and the relationship of Slepian matrix and Dirichlet kernel, we obtain the claimed result. 

This theorem echoes the numerical findings in Figure \ref{fig:Vplot}. We shall mention that this result is critical for the detection resolution of Thomson's multitaper $F$-test, which will be reported in the other paper.

\section{Mathematical Background}\label{section Mathematical Background}

\subsection{A quick review of DPSWF and its relationship with DPSS tapers}\label{app_tapers}
We quickly review DPSS and  DPSWF and summarize their basic properties.  
More properties can be found in \cite{slepian1978pswfV,waldenspectral}. 

Recall that for a discrete sequence of real numbers $g:=\{g_j\}_{j\in \ZZ}$ such that $\sum_{j=-\infty}^{\infty} g_j^2 = 1$, its Fourier transform is defined as $G(f) = \sum_{-\infty}^{\infty} g_je^{-i2\pi jf}$. By Parseval's Theorem we have $\int_{-\frac{1}{2}}^{\frac{1}{2}}|G(f)|^2df=1$. 
It is well-known that a simultaneously band- and time-limited nontrivial signal cannot exist due to the Paley-Wiener theorem. 
A natural question to ask is if we could find a function that is as close as possible to such simultaneously band- and time-limited signal. This can be quantified by two types of concentration. For the first type, among all band-limited $g$ such that $G(f)=0$ for $|f|>W$, where $W\in (0,1/2)$, find one that maximizes $\sum_{j=0}^{N-1}g_j^2$ for a given $N\in \mathbb{N}$; that is, we need solve the eigenvalue problem of the following integral equation:
\begin{align}
    \int_{-W}^{W}\frac{\sin(N\pi(f'-f))}{\sin(\pi(f'-f))}G(f)df=\lambda G(f')\text{ for }|f'|\leq W\label{band_eigen}\,.
\end{align}
As discussed in \cite[(103a)]{waldenspectral}, there are $N$ eigenfunctions to \eqref{band_eigen}, called the DPSWF, denoted by $U_0(f;N,W),\ldots,U_{N-1}(f;N,W)$. The eigenvalues are ordered by 
\[
1>\lambda_0(N,W)>\ldots>\lambda_{N-1}(N,W)>0\,. 
\]
We extend $U_k(f;N,W)$ to $[-\frac{1}{2},\frac{1}{2}]$ using \eqref{band_eigen} and normalize them such that 
\[
\int_{-\frac{1}{2}}^{\frac{1}{2}}U_k(f;N,W)U_l(f;N,W)=\delta_{kl},\ U_k(0;N,W)\geq 0,\ \mbox{ and }\ \frac{dU_k(f;N,W)}{df}\Big|_0\geq 0\,. 
\]
The DPSWFs also satisfy 
\begin{equation}
\int_{-W}^{W}U_k(f;N,W)U_l(f;N,W)df=\delta_{kl}\lambda_l(N,W)\,. \label{DPSWF behavior in the frequency domain}
\end{equation}
We call $W$ the {\em bandwidth}, $[-W,W]$ the {\em inside band} (or passband), $[-1/2,-W)\cup(W,1/2]$ {\em outside band}, and any small spectral range including $W$ a {\em critical band}. 
The $j$-th DPSS, denoted as $v_j=v_j(N,W)$, then is defined as \cite[page 104]{waldenspectral}
\begin{equation}
v_{j,k}(N,W) = \frac{(-1)^j}{\epsilon_j\lambda_k(N,W)}\int_{-W}^W U_j(f;N,W)e^{i2\pi f(k-(N-1)/2]}df\,,\label{definition jth DPSS in l^2}
\end{equation}
where $k\in \mathbb{Z}$ and $\epsilon_j=1$ for even $j$ and $\epsilon_j=\sqrt{-1}$ for odd $j$. All $v_j$ are real-valued \cite[page 104]{waldenspectral}. 
We also have \cite[(3.12) on page 122]{waldenspectral}
\begin{align}
U_j(f;N,W)&\,=(-1)^j\epsilon_j\sum_{k=0}^{N-1}v_{j,k}(N,W)e^{-i2\pi f[k-(N-1)/2]}\label{definition jth DPWSF from Slepian} \\
&\,=\epsilon_j\sum_{k=0}^{N-1}v_{j,k}(N,W)e^{i2\pi f[k-(N-1)/2]}\nonumber
\end{align}
for $f\in [-1/2,1/2]$.
From now on, unless necessary, we avoid $N,W$ in the notation for simplicity.

The second way to generate the $k$-th order DPSS is solving the following eigenvalue problem \cite[page 104]{waldenspectral}:
\begin{align}
    \mathbf A(g_0,\ldots,g_{N-1})^\top = \lambda (g_0,\ldots,g_{N-1})^\top\,,\label{time_eigen}
\end{align}
where $\mathbf A\in \mathbb{R}^{N\times N}$ is the {\em prolate matrix} defined as
\begin{align}\label{definition of A matrix}
\mathbf A_{i+1,j+1}:=\frac{\sin(2\pi W(i-j))}{\pi(i-j)}\text{ for }0\leq i,j\leq N-1\,.
\end{align}
The $j$th eigenvalue of $\mathbf A$ is $\lambda_j$, and the associated $j$th eigenvector is denoted as $\vv_j\in \mathbb{R}^N$. 
$\vv_j$ is a subsequence of $v_j$ with indices in the range $0$ to $N-1$ \eqref{definition jth DPSS in l^2} \cite[page 104]{waldenspectral}. We call $\vv_j$ the $j$th {\em Slepian taper}, or simply the $j$-th DPSS since $v_j$ can be recovered from $\vv_j$ \cite[page 104]{waldenspectral} (or see \cite{slepian1978pswfV}).
For the Slepian tapers, we follow the following polarity convention \cite{waldenspectral}. 
For symmetric tapers ($j=0,2,4\ldots$), assume
\[
\sum_{k=0}^{N-1}\vv_{j,k}>0\,,
\]
and for skew-symmetric tapers ($k=1,3,5\ldots$), the taper is made to start with a positive first lobe; that is,
\[
\sum_{k=0}^{N-1}(N-1-2k)\vv_{j,k}>0\,.
\]

Note that DPSS' are the eigenvectors of $\mathbf A$, which is a symmetric Toeplitz matrix. Because $\mathbf A$ and the domain are symmetric, if a sequence is a solution, its time-reversed version is also a solution. Consequently, the eigenfunctions partition into purely symmetric and purely skew-symmetric functions.
By standard ordering conventions, the DPSS alternate in symmetry so that even order DPSS tapers are even, and odd order DPSS tapers are odd. Clearly, $\sum_{j=0}^{N-1}\vv_{k,j}=0$ when $k$ is odd.

$\mathbf A$ commutes with the {\em Slepian's commuting matrix} $\mathbf D=\mathbf D(W)\in \mathbb{R}^{N\times N}$ \cite[(14)]{slepian1978pswfV}, which is symmetric tridiagonal with nonzero entries 
\[
\mathbf D_{l,l}=\Big(\frac{N-1}{2}-l\Big)^2\cos(2\pi W)\,,\ \mathbf D_{l,l-1} =\frac{l(N-l)}{2}\,,\ \mbox{ and }\mathbf D_{l,l+1}=\frac{(l+1)(N-1-l)}{2}\,.
\]
Since $\mathbf A$ and $\mathbf D$ commute, they share the same eigenvectors, Slepian tapers.  Denote its eigenvalue decomposition as $\mathbf D\vv_k=\theta_k\vv_k$, where $\theta_k$ is the eigenvalue ordered decreasingly. The tridiagonal eigenvalue problem can be written row by row as a difference equation \cite[(25)]{slepian1978pswfV}
\[
\frac{n(N-n)}{2}\vv_{k,n-1}+\Big(\frac{N-1}{2}-n\Big)^2\cos(2\pi W)\vv_{k,n}+\frac{(n+1)(N-1-n)}{2}\vv_{k,n+1}=\theta_k\vv_{k,n}\,.
\]
Denote $c_n:=n-(N-1)/2$ as the phase shift in the shifted Fourier transform. We have
\[
\frac{(n+1)(N-1-n)}{2}=-\frac{c_n(c_n+1)}{2}+\frac{N^2-1}{8}
\]
and
\[
\frac{n(N-n)}{2}=-\frac{c_n(c_n-1)}{2}+\frac{N^2-1}{8}\,.
\]
By plugging these two equations into the difference equation, collecting terms with $c_n^2$, $c_n$, and others, and observing that their DTFTs satisfy
\[
\sum_{n=0}^{N-1}\Big[-\frac{1}{2}c_n^2(\vv_{k,n-1}+\vv_{k,n+1})+\cos(2\pi W)c_n^2 \vv_{k,n}\Big] e^{-i2\pi c_n f}=(\cos(2\pi f)-\cos(2\pi W))U_k''(f)\,,
\]
\[
\sum_{n=0}^{N-1}\Big[\frac{1}{2}c_n(\vv_{k,n-1}-\vv_{k,n+1})\Big]e^{-i2\pi c_n f}=-\sin(2\pi f) U_k'(f)\,,
\]
and
\[
\sum_{n=0}^{N-1}\Big[\frac{N^2-1}{8} (\vv_{k,n-1}+\vv_{k,n+1})\Big]e^{-i2\pi c_n f}=\Big(\frac{N^2-1}{4}\Big)\cos(2\pi f) U_k(f)\,.
\]
we see that the DPSWF satisfies the following differential equation over $f \in (-1/2, 1/2)$ \cite[(113)]{slepian1978pswfV}:
\begin{align}
\frac{d}{df}[p(f) U_k'(f)] + 4\pi^2 Mq_k(f)U_k(f)=0\,,\label{2nd order ODE for DSPWF}
\end{align}
where 
\[
p(f):=\cos(2\pi f) - A,\, q_k(f):=\cos(2\pi f) - B_k\,,
\]
\[
A := \cos(2\pi W)\,,\ B_k := \theta_k/M\in (-1,1)\,, \ \mbox{ and }\ M := (N^2-1)/4\,.
\] 
Geometrically, $B_k$ quantifies the {\em turning point} of the DPSWF. It describes the boundary where the wave function transits from oscillation to its exponential decay to $0$. 
This ODE has an interesting structure. When $k<2NW$, over $[0,W]$, 
$p(f)>0$ and has simple zero of the leading coefficient and $q_k$ has a turning point at $w_k$, 
and over $[W,1/2]$, $p(f)<0$ when $f\in (W,1/2]$ and $p(W)=0$, and $q(f)<0$, so that $U_k$ oscillates again with a simple zero of the leading coefficient.

In this paper, we call $k\leq 2NW-C(N)$ for some $C(N)>0$ satisfying $\frac{C(N)}{2NW}\to 0$ as $N\to \infty$ the {\em pre-plunge} region, 
$2NW-C(N)<k<2NW+C(N)$ the {\em plunge (or transition)} region, and $k\geq 2NW+C(N)$ the {\em post-plunge} region.

\subsection{Existing results regarding DPSWF}

The following theorem is a simplified version of the result in \cite[Section 3.3]{karnik2020improvedchardpsseigenval}, tailored to our purposes. Note that, as $W\to 0$ with $N\to \infty$, we have$100NW+25<4N$ for sufficiently large $N$, which justifies this simplification.   

\begin{theorem}\label{karnik_bound}
Assume when $W\to 0$ as $N\to \infty$. When $N$ is sufficiently large and $k\leq \lfloor 2NW\rfloor -1$, we have 
\[
\lambda_k \geq 1-10 \exp\Big(-\frac{\lfloor 2NW\rfloor-k-7}{\frac{2}{\pi^2}\log(100NW+25)}\Big)\,,
\]
and when $k\geq \lfloor 2NW\rfloor$, we have 
\[
\lambda_k \leq 10 \exp\Big(-\frac{k-\lfloor 2NW\rfloor-6}{\frac{2}{\pi^2}\log(100NW+25)}\Big)\,.
\]
For $\epsilon\in (0,1/2)$, when $N$ is sufficiently large, we have \cite[Theorem 2]{karnik2020improvedchardpsseigenval}
\[
\#\{k|\, \epsilon\leq \lambda_k\leq 1-\epsilon\}\leq \frac{2}{\pi^2}\log(100NW+25)\log\Big(\frac{5}{\epsilon(1-\epsilon)}\Big)+7\,.
\] 
\end{theorem}

Note that by taking the trace of the matrix $A$ in \eqref{time_eigen}, we have $\sum_{k=0}^{N-1}\lambda_k = 2NW$. 
The following lemma is a direct generalization of \cite[Corollary 2]{karnik2020improvedchardpsseigenval} that we will use later.

\begin{lemma}\label{karnik_bound22}
Assume $W\to 0$ when $N\to \infty$. For any $p> 0$, 
\begin{align}
\sum_{k=0}^{K-1}(1-\lambda_k)^p &\,\leq \frac{20}{p\pi^2}\log(100NW+25)\exp \Big(-\frac{\lfloor 2NW\rfloor-K-7}{\frac{2}{p\pi^2}\log(100NW+25)}\Big)\,.\label{bound sum of sqrt root of 1-lambdak}
\end{align}
Take $K=\lfloor 2NW-\frac{2c}{p\pi^2}\log(100NW+25)\log(N)\rfloor$ for some $c>0$. For any $c'<c$, we have $\sum_{k=1}^K(1-\lambda_k)^p\leq N^{-c'}$ when $N$ is sufficiently large. Similar results hold for $\sum_{k=K}^{N-1}\lambda_k^p$ for $K>2NW$.
\end{lemma}

\begin{proof}

The bound comes from exactly the same proof in \cite[Section 5.3]{karnik2020improvedchardpsseigenval} by modifying constants.
When $N$ is sufficiently large, we have the control $\frac{20}{p\pi^2}\log(100NW+25)N^{-c}\leq N^{-c'}$ for $c'<c$.

\end{proof}

The following theorem summarized the results from \cite{said2023non} that are relevant to our analysis. We refer interested readers to \cite{said2023non} for additional results. Note that the ODE for DPSWFs in \cite{said2023non} is defined on $(-\pi,\pi)$, whereas our analysis focuses on $(-1/2,1/2)$. Accordingly, the results below are presented in a rescaled form of those in \cite{said2023non}.

\begin{theorem}\label{theorem:said2023non:summary}
For DPSWF, denote $\Omega_i$ to be zeros of $U_k$ and $\nu_i$ to be zeros of $U_k'$, and set $\Omega_0=0$, $\Omega_{0.5}=0.5$, $\nu_0=0$, and $\nu_{0.5}=0.5$. Denote $\nu_k^*$ to be the maximal extremal point of $U'_k$ if $\frac{1}{2\pi}\cos^{-1}(B_k)<W$, and $\nu_k^*=W$ if $\frac{1}{2\pi}\cos^{-1}(B_k)>W$ \cite[(9)]{said2023non}. We have
\begin{itemize}
\item[(a)] When $k\leq 2NW-1$, $B_k>A$; when $2NW\leq k\leq N-1$, $B_k<A$; when $2NW-1<k<2NW$, either equality is possible. \cite[Theorem 1]{said2023non}

\item[(b)] For any two arbitrary extremum points $x$ and $y$ of $U_k$ in $(-W, W)$, if $|x|<|y|$, then $|U_k(x)|<|U_k(y)|$. If further $B_k<A$, then $|U_k(x)|<|U_k(y)|<|U_k(W)|$. \cite[Theorem 2]{said2023non}

\item[(c)] 
If $k\leq 2NW-1$, then \cite[Theorem 2]{said2023non}
\[
\max_{f\in[0,W]}|U_k(f)|=|U_k(\nu_k^*)|\ \mbox{ and }\ \max_{f\in [W,1/2]}|U_k(f)|=|U_k(W)|\,.
\] 
If $k\geq 2NW$, $\sup_{f\in [0,W]}|U_k(f)|=|U_k(W)|$ \cite[Corollary 2]{said2023non}.

\item[(d)] For $f\in \big(0,\frac{1}{2\pi}\min(\cos^{-1}(A),\cos^{-1}(B_k)\big)$, we have \cite[Theorem 3 and (54)]{said2023non}
\begin{align*}
\sqrt{\frac{\cos(2\pi f)-A}{1-A}\frac{\cos(2\pi f)-B_k}{1-B_k}} |U_k(f)|^2\leq  |U_k(0)|^2+\frac{1-A}{1-B_k}\frac{1}{4\pi^2M}|U_k'(0)|^2\,.
\end{align*}

\item[(e)] We have \cite[Theorem 4]{said2023non}
\[
\sup_{f\in [-1/2,1/2]}|U_k(f)|\leq 
\left\{
\begin{array}{ll}
\mathsf U_k(A,B_k,0)& \mbox{ when } B_k>A\\
\mathsf U_k(-A,-B_k,0)& \mbox{ when } B_k<A
\end{array}
\right.
\]
where $\mathsf U_k(x,y,a):=\min_{i\in\{1,2\}}\kappa_i(x,y)\sqrt{|U_k(a)|^2+\frac{1-x}{1-y}\frac{1}{4\pi^2M}|U_k'(a)|^2}(\cos^{-1}(x)-\cos^{-1}(y))^{1/4}$,
$\kappa_1(x,y):=\frac{9}{8}\big(\frac{\sqrt{9M\frac{1-y}{1-x}}}{\cos^{-1}(x)-\cos^{-1}(y)}\big)^{1/4}$, and $\kappa_2(x,y):=\frac{5}{4}\big(5M\frac{1-y}{1-x}\big)^{1/4}$.

\item[(f)] $U_k$ has $N-1$ roots $\Omega_i$ in $(-1/2,1/2]$ and $U'_k$ has $N-1$ roots $\nu_i$ interlaced with the zeros of $U_k$ \cite[(8)]{said2023non}. When $B_k>A$ and $k$ is odd \cite[Lemma 2 and (7)]{said2023non},
\begin{align*}
\Omega_0<&\,|\nu_{\pm1}|<|\Omega_{\pm1}|<\ldots<\underbrace{|\nu_{\pm (k+1)/2}|<\frac{1}{2\pi}\cos^{-1}(B_k)<W}_{I_o}\\
&\,<|\Omega_{\pm (k+1)/2}|<\ldots\left\{\begin{array}{ll}
|\nu_{\pm \frac{N-1}{2}}|<\Omega_{0.5}& \mbox{for } N \mbox{ odd}\\
|\Omega_{\pm \frac{N-2}{2}}|<\nu_{0.5}& \mbox{for } N \mbox{ even}
\end{array}\right.
\end{align*}
and when $B_k>A$ and $k$ is even,
\begin{align*}
\nu_0<&\,|\Omega_{\pm1}|<|\nu_{\pm1}|<\ldots<\underbrace{|\nu_{\pm k/2}|<\frac{1}{2\pi}\cos^{-1}(B_k)<W}_{I_e}\\
&\,<|\Omega_{\pm k/2+1}|<\ldots\left\{\begin{array}{ll}
|\nu_{\pm \frac{N-2}{2}}|<\Omega_{0.5}& \mbox{for } N \mbox{ even}\\
|\Omega_{\pm \frac{N-1}{2}}|<\nu_{0.5}& \mbox{for } N \mbox{ odd}
\end{array}\right.\,.
\end{align*}
When $B_k<A$, $I_o$ is replaced by $W<\frac{1}{2\pi}\cos^{-1}(B_k)<|\nu_{\pm (k+1)/2}|$, and $I_e$ is replaced by $W<\frac{1}{2\pi}\cos^{-1}(B_k)<|\nu_{\pm k/2}|$. In other words, $U_k$ has $k$ roots in $(-W, W)$.
\end{itemize}
\end{theorem}

This theorem provides detailed information about DPSWF's morphology. For example, when $k<2NW$, $|U_k(\nu_k^*)|=\max_{f\in(0,W)}|U_k(f)|$ and there is no extreme points or zeros between $\nu_k^*$ and $W$, or $\frac{1}{2\pi}\cos^{-1}(B_k)$ is not an extreme point nor root, and the magnitude of $U_k$ is decreasing from $\nu_k^*$ to $W$. 
Since $\nu_k^*<\frac{1}{2\pi}\cos^{-1}(B_k)$, $\cos(2\pi \nu_k^*)>B_k>A$.

Theorem \ref{theorem:said2023non:summary}(a) identifies a potentially undetermined comparison between $B_k$ and $A$ when $2NW$ is not an integer, namely $k=\lfloor 2NW\rfloor$. Since there is at most one such index, and the proof is already lengthy, we impose the following assumption for the remainder of the analysis. The other case can be handled similarly and is omitted for brevity.

\begin{Assumption}
When $2NW-1<k<2NW$, assume $B_k<A$.
\end{Assumption}

With this assumption, in the following paper, we discuss $k<2NW$ and $k\geq 2NW$.

The main challenge of applying results from  \cite{said2023non}, particularly Theorem \ref{theorem:said2023non:summary}(d), to control $U_k$ size is that it blows up when $f$ is close to $w_k$ or $W$. When $k\geq 2NW$ and $k$ is even, we have $B_k<A$ by \ref{theorem:said2023non:summary}(a), so we can apply \ref{theorem:said2023non:summary}(d) for $f\in (0,W)$ and have
$ |U_k(f)|^2\leq  |U_k(0)|^2\sqrt{\frac{1-A}{\cos(2\pi f)-A}\frac{1-B_k}{\cos(2\pi f)-B_k}}$.
Since $\cos(2\pi f)-A$ is a monotonically decreasing function, for a large but fixed constant $c_*$, consider $f_*=\frac{1}{2\pi}\cos^{-1}(c_*^{-1}+(c_*-1)c_*^{-1}A)$ so that $\frac{1-A}{\cos(2\pi f)-A}\leq c_*$. Also, use the property that $B_k<A$ implies $\frac{1-B_k}{\cos(2\pi f)-B_k}<\frac{1-A}{\cos(2\pi f)-A}$, we obtain
\begin{align}
 |U_k(f)|^2\leq  c_*|U_k(0)|^2=O(\lambda_kW^{-1})\,,\label{application of Said 2023 control for inside band for k>2NW}
\end{align}
which captures the fact that $U_k$ is almost like $0$ on the inside band, but it only captures this over $[0,f_*]$.
Note that when $c_*=10$, $f_*=0.9486W$. However, this bound blows up when $f$ approaches $W$. Moreover, it has limited applicability when $k<2NW$, since in this case $B_k>A$, restricting the spectral range over which the bound is useful. For example, when $k$ is small, $B_k$ is almost $1$, so the bound can only control a spectral range near $0$ that is much narrower than $[0,W]$.

Similarly, we can capture the fact that, when $k<2NW$, $U_k$ is nearly $0$ on the outside band. To see this, we apply the flipping trick \cite[(13)]{slepian1978pswfV}.
Denote $U_k(f;N,W)$ and $\lambda_k(N,W)$ to be the $k$-th DPSWF with number $N$ and bandwidth $W$ and the associated eigenvalue. We know $U_k(1/2-f;N,W)$ is the $k'$-th DPSWF with length $N$ and parameter $W':=1/2-W$, where $k':=N-k-1$; that is, $U_k(1/2-f;N,W)=U_{k'}(f;N,W')$. Clearly, $k'=N-k-1>N-2NW= 2NW'$ (post-plunge region). In this case, $-B_k=B_k'<A'=\cos(2\pi W')=-A$, so for $f\in [0,W')$, we have
\begin{align*}
&|U_k(1/2-f;N,W)|^2=
 |U_{k'}(f;N,W')|^2 \\
 \leq \,& \sqrt{\frac{1-A'}{\cos(2\pi f)-A'}\frac{1-B'_k}{\cos(2\pi f)-B'_k}}\Big(|U_{k'}(0;N,W')|^2+\frac{1-A'}{1-B'_k}\frac{1}{4\pi^2M}|U_{k'}'(0;N,W')|^2 \Big) \\
 =\,&  \sqrt{\frac{1+A}{\cos(2\pi f)+A}\frac{1+B_k}{\cos(2\pi f)+B_k}}\Big(\Big|U_{k}\Big(\frac{1}{2};N,W\Big)\Big|^2+\frac{1+A}{1+B_k}\frac{1}{4\pi^2M}\Big|U_{k}'\Big(\frac{1}{2};N,W\Big)\Big|^2 \Big)\,.
\end{align*}
As a result, since $1>B_k>A=\cos(2\pi W)$, when $f\in (W,1/2]$, 
$|U_k(f;N,W)|^2 \leq  \sqrt{\frac{1+A}{A-\cos(2\pi f)}\frac{1+B_k}{B_k-\cos(2\pi f)}}\big(\big|U_{k}\big(\frac{1}{2};N,W\big)\big|^2+ \frac{1}{4\pi^2M}\big|U_{k}'\Big(\frac{1}{2};N,W\big)\big|^2 \big)$, 
Note that $\frac{1+A}{A-\cos(2\pi f)}$ dominates $\frac{1+B_k}{B_k-\cos(2\pi f)}$.
Since $A-\cos(2\pi f)$ is an increasing function, for $c_*>0$ a large constant, we can find $f_*$ such that $1+A=c_*(A-\cos(2\pi f))$, so that when $f\in \big(\frac{\cos^{-1}(\frac{c_*-1}{c_*}A-\frac{1}{c_*})}{2\pi}, 1/2\big]$, 
\begin{align}
|U_k(f;N,W)|^2 \leq   c_*\Big(\Big|U_{k}\Big(\frac{1}{2};N,W\Big)\Big|^2+ \frac{1}{4\pi^2M}\Big|U_{k}'\Big(\frac{1}{2};N,W\Big)\Big|^2 \Big)\,.\label{application of Said 2023 control for outside band for k<2NW}
\end{align}
Here, when $c_*=20$, we have $\frac{\cos^{-1}(\frac{c_*-1}{c_*}A-\frac{1}{c_*})}{2\pi}\approx 3.492 W$. 
Therefore, to get the claim that $U_k$ is nearly $0$ over $f\in \big(\frac{\cos^{-1}(\frac{c_*-1}{c_*}A-\frac{1}{c_*})}{2\pi}, 1/2\big]\subset [0,1/2]$ when $k<2NW$, depending on the sign of $N$, we need a control of $U_k$ or $U'_k$ at $1/2$.

To obtain a bound that remains effective near $W$, we may consider the global bound in Theorem \ref{theorem:said2023non:summary}(e). When $k<2NW$, we have $\kappa_1(A,B_k)=\frac{9}{8}\big(\sqrt{9M\frac{1-B_k}{1-A}}\frac{1}{2\pi (W-w_k)}\big)^{1/4}$ and $\kappa_2(A,B_k)=\frac{5}{4}\big(4M\frac{1-B_k}{1-A} \big)^{1/4}$. Assume that $\kappa_1$ is minimal, which holds when $W\to 0$ and $NW\to \infty$ when $N\to \infty$ since $\frac{\kappa_1(A,B_k)}{\kappa_2(A,B_k)}=O(\frac{1}{NW})$ when $N\to\infty$. Since we only need to consider even $k<2NW$, it suffices to control the size of $\big( \frac{(1-B_k)N^2}{1-A} \big)^{1/8}$. Although this bound is not ideal for obtaining optimal pointwise control, it is nevertheless useful for our purposes.

\section{Technical lemmas}\label{section technical lemmas}

We start with preparing some lemmas to control the relationship between $B_k$ and $A$ when $k<2NW$. 

\begin{lemma}\label{lemma: control M(1-Bk) asymptotically}
For any $W\in (0,1/2)$ and $N\geq 1$, we have
\begin{equation}
    \cos^2(\pi W) \chi_k(\tilde{c}) \leq 2M(1-B_k) \leq \chi_k(\tilde{c})\,,\label{proof: asymptotic control of ev of PSWF and DPSWF}
\end{equation}
where $\tilde{c}^2 := \frac{4M\sin^2(\pi W)}{\cos(\pi W)}$ is the {\em continuous time-bandwidth parameter} and $\chi_k(\tilde c)$  is the $k$-th eigenvalue of the PSWF ODE with the time-bandwidth parameter $\tilde{c}$.
\end{lemma}

\begin{proof}

Recall that the ODE of the continuous PSWF with time-bandwidth parameter $\tilde{c}$ is \cite{osipov2013prolate}
\[
-((1-x^2)\psi_k'(x))'+c^2x^2\psi_k(x)=\chi_k(\tilde c)\psi_k(x)\,.
\] 
The eigenvalue satisfies
\[
\chi_k(\tilde c)=\min_{S_k}\max_{\psi\in S_k\backslash\{0\}}R_0(\psi)\,,  \mbox{ where } R_0(\psi):=\frac{\int_{-1}^1[(1-x^2)\psi'(x)^2+{\tilde c}^2x^2\psi(x)^2]dx}{\int_{-1}^1\psi(x)^2dx}
\]
and $S_k$ is a $k$-dim subspace of $L^2[-1,1]$.
For the ODE for DPSWF \eqref{2nd order ODE for DSPWF},
consider a change of variable $x = \frac{\sin(\pi f)}{\sin(\pi W)}$. 
As $f$ varies in $[-W, W]$, $x$ maps exactly to $[-1, 1]$. Using the identities $\cos(2\pi f) =1-2\sin^2(\pi f)= 1 - 2x^2\sin^2(\pi W)$ and $A = 1 - 2\sin^2(\pi W)$, the ODE becomes
\[
-(\tilde p(x)\psi'_k(x))'+\tilde q(x)\psi_k(x)=\tilde\lambda_kw(x)\psi_k(x)
\] 
on $[-1, 1]$ with 
\begin{align*}
\tilde p(x):=\,&(1-x^2)\sqrt{1-x^2\sin^2(\pi W)}\,,\\ 
\tilde q(x):=\,&\frac{4M\sin^2(\pi W)x^2}{\sqrt{1-x^2\sin^2(\pi W)}}=\tilde{c}^2x^2\frac{\cos(\pi W)}{\sqrt{1-x^2\sin^2(\pi W)}}\\
\tilde w(x):=\,&\frac{1}{\sqrt{1-x^2\sin^2(\pi W)}}\,,\ \tilde\lambda_k:=2M(1-B_k)\,.
\end{align*}
The eigenvalue $\tilde\lambda_k$ satisfies 
\[
\tilde\lambda_k=\min_{S_k}\max_{\phi\in S_k\backslash\{0\}}R_1(\phi)\,, \mbox{ where } R_1(\phi):=\frac{\int_{-1}^1[\tilde p(x)\phi'(x)^2+\tilde q(x)\phi(x)^2]dx}{\int_{-1}^1\tilde w(x)\phi(x)^2dx}\,,
\]
The following bounds are satisfied:
\begin{align*}
(1-x^2)\cos(\pi W)\leq \tilde{p}(x)\leq 1-x^2\,,\ 
\tilde{c}^2x^2\cos(\pi W)\leq \tilde{q}(x)\leq \tilde{c}^2x^2\,, \mbox{ and }
1\leq \tilde{w}(x)\leq \frac{1}{\cos(\pi W)}\,,
\end{align*} 
so that 
\[
\int_{-1}^1[\tilde p(x)\phi'(x)^2+\tilde q(x)\phi(x)^2]dx\leq \int_{-1}^1[(1-x^2)\phi'(x)^2+{\tilde c}^2x^2\phi(x)^2]dx
\]
and
\[
\int_{-1}^1\tilde w(x)\phi(x)^2dx\geq \int_{-1}^1\phi(x)^2dx\,,
\]
which gives $R_1(\phi)\leq R_0(\phi)$ for every $\phi\neq 0$. Similarly, we have $R_1(\phi)\geq \cos^2(\pi W)R_0(\phi)$. 
By the Courant min-max theorem, we obtain the claimed bound
\[
\cos^2(\pi W) \chi_k(\tilde{c}) \leq 2M(1-B_k) \leq \chi_k(\tilde{c})\,.
\]
\end{proof}

We prepare some basic notation and facts for the next lemma. Denote 
\[
K(x):=\int_0^{\pi/2}\frac{1}{\sqrt{1-x^2\sin^2(\phi)}}d\phi\ \mbox{ and } \ E(x):=\int_0^{\pi/2}\sqrt{1-x^2\sin^2(\phi)}d\phi\,, 
\]
where $x\in [0,1]$, to be complete elliptic integrals of the first and second kinds. $E$ is strictly decreasing with $E(0)=\pi/2$ and $E(1)=1$. $K$ is strictly increasing, blows up when $x\to 1$, with $K(0)=\frac{\pi}{2}$ and $K(x)\sim \log\big( \frac{4}{\sqrt{1-x^2}} \big)$ when $x\to 1$ so that $(1-x^2)K(x)\to 0$ when $x\to 1$. Denote 
\[
H(x):=E(x)-(1-x^2)K(x)\,, 
\]
$x\in [0,1]$. We have $H(0)=0$ and $H(1)=1$. Recall that 
\[
E'(t)=\frac{E(t)-K(t)}{t}\ \mbox{ and }\ K'(t)=\frac{E(t)}{t(1-t^2)} -\frac{K(t)}{t}
\] 
for $t\in [0,1)$. We therefore have $H'(t)=tK(t)>0$ and $H''(t)=\frac{E(t)}{1-t^2}>0$, so $H(x)$ is strictly increasing and convex.

\begin{lemma}\label{taper_upper bound in frequency domain, control of 1-A/B-A}
Assume $W\to 0$ and $NW\to\infty$ when $N\to \infty$. 
When $2\leq k< 2NW-1$ and $N$ is sufficiently large, we have the following controls
\begin{enumerate}
\item $\displaystyle 1-A= 2\pi^2W^2(1+O(W^2))$ and $\displaystyle W= \frac{\sqrt{1-A}}{\sqrt{2}\pi} (1+O(W^2))$;
\item $\displaystyle
\frac{\pi^2 W(2NW-\frac{4}{\pi}(k+1))}{N} \leq B_k-A\leq \frac{\pi^2 W(2NW-k)}{N-1}$
for all $k<2NW$, and when $\frac{3}{2}NW<k<2NW-1$, we have 
\begin{align*}
&2(\pi W)^2
\Big( 
0.694 \Big(\frac{2NW}{k+1} \Big)^2 - 1.08 \Big(\frac{2NW}{k+1} \Big) + 0.386
\Big) < B_k-A\\
< \,& 2(\pi W)^2 \Big( 
-0.155 \Big(\frac{2NW}{k} \Big)^2 +0.901 \Big(\frac{2NW}{k} \Big) - 0.746
\Big)\,;\nonumber
\end{align*}

\item $\displaystyle
 \frac{\pi^2 kW}{N-1}\leq 1-B_k\leq \frac{4\pi (k+1)W}{N}$
for all $k<2NW$, and when $\frac{3}{2}NW<k<2NW-1$, we have
\begin{align*}
&2(\pi W)^2
 \Big( 
0.155 \Big(\frac{2NW}{k} \Big)^2 -0.901 \Big(\frac{2NW}{k} \Big) + 1.746
\Big)
< 1-B_k\\
< \,& 2(\pi W)^2 \Big( 
-0.694 \Big(\frac{2NW}{k+1} \Big)^2 +1.08 \Big(\frac{2NW}{k+1} \Big) + 0.614
\Big)\,.\nonumber
\end{align*}

\end{enumerate}

\end{lemma}

This Lemma tells us that when $k\sim 1$, $1-B_k\sim \frac{W}{N}$, and $B_k-A\sim W^2$. On the other hand, when $k\sim 2NW(1-W^\theta)$ for $\theta>0$, $1-B_k\sim W^2$, and $B_k-A\sim W^{2+\theta}$. Note that if $W\to 0$ when $N\to \infty$, we have
\begin{align}
 w_k^2= \frac{1-B_k}{2\pi^2} (1+O(w_k^2))\ \mbox{ and }\   W^2-w_k^2=\frac{B_k-A}{2\pi^2}(1+O(W^2))\,.\label{proof of Ukf decay, wk and W in Bk and A}
\end{align}
Moreover, combined with Theorem \ref{theorem:said2023non:summary}(e), we have
\begin{align}\label{Uk worst case bound}
\sup_{f\in [-1/2,1/2]}|U_k(f)|\leq \frac{3^{9/4}}{2^{11/2}\pi^{1/8}}\frac{(k+1)^{1/8}N^{1/8}}{W^{5/8}}\,.
\end{align}

\begin{proof}
Denote $c=\pi NW$. By the definition of $\tilde{c}$ in Lemma \ref{lemma: control M(1-Bk) asymptotically}, we have $\tilde{c}/c =1+O(W^2)\to 1$ as $N\to \infty$. 
Since $k<2NW$, we have \cite[Theorem 3.1]{osipov2013prolate}
\[
\chi_k(c)<c^2\,.
\]
Since $\chi_k(c)$ is continuous \cite[(3.3)]{osipov2013prolate}, $\chi_k(\tilde{c}) < \tilde{c}^2$ when $N$ is sufficiently large. Clearly,  $\tilde{c}/M\to 0$ when $N\to \infty$.
Next, from \eqref{proof: asymptotic control of ev of PSWF and DPSWF}, 
when $N$ is sufficiently large,
\[
1-B_k=\frac{\chi_k(c)}{2M}(1+O(W^2))\,.
\]
Similarly, by the definition of $\tilde{c}^2$, since $\frac{2M(1-A)}{\tilde c^2} = \frac{4M\sin^2(\pi W)}{\tilde c^2}=\cos(\pi W)=1+O(W^2)$, when $N$ is sufficiently large, we have 
\[
1-A=\frac{c^2}{2M}(1+O(W^2))\,.
\]
Hence,
\begin{equation}
    \frac{1-A}{1-B_k} = \frac{c^2}{\chi_k(c)}(1+O(W^2)) \ \mbox{ and }\ 
    B_k-A= \frac{c^2-\chi_k(c)}{2M}(1+O(W^2))\,.
\end{equation}
By \cite[Proposition 2]{bonami2016uniform}, for all $2\leq k<2NW-1$ we have 
\begin{align}
c(c-2k-2)\leq c^2-\chi_k(c)\leq c(c-\pi k/2)\,, \label{bound of c2-chik immediately from bonami prop2}
\end{align}
and hence
\[
\frac{\pi^2 W(2NW-\frac{4}{\pi}(k+1))}{N}(1+O(W^2)) \leq B_k-A\leq \frac{\pi^2 W(2NW-k)}{N-1}(1+O(W^2))\,.
\] 
However, the lower bound becomes trivial when $k> \frac{\pi}{2} NW-1$. This is because the bound in \cite[Proposition 2]{bonami2016uniform} is global. For our purpose, we need a finer lower bound when $k>\frac{3}{2} NW$. 
Below we establish a sharper bound when $k> \frac{3}{2} NW$.

Denote
\[
Q(\kappa):=\int_0^{\min\{1,1/\kappa\}} \sqrt{ \frac{1-\kappa^2t^2}{1-t^2} } dt
\]
is a decreasing function on $[0,\infty)$.
Recall $\Phi$ above the statement of \cite[Theorem 1]{bonami2016uniform}, which is the inverse function of the monotonically increasing function $\kappa\mapsto \frac{\kappa}{Q(\kappa)}$. In other words, 
\[
\Phi\Big( \frac{\kappa}{Q(\kappa)} \Big)=\kappa=\Phi(s)
\] 
if we denote $s:= \frac{\kappa}{Q(\kappa)}$, which is an increasing function on $[0,\infty)$. Therefore, $\Phi$ is a monotonically increasing function with $\Phi(0)=0$ and $\Phi(1)=1$. 
Note that when $\kappa>1$, with change of variable $\kappa t=\sin(\phi)$, we have 
\[
Q(\kappa)=\frac{1}{\kappa}\int_0^{\pi/2}\frac{\cos^2(\phi)}{\sqrt{1-\kappa^{-2}\sin^2(\phi)}}d\phi\,.
\] 
Therefore, $\kappa Q(\kappa)$ is decreasing on $[1,\infty)$. By a direct change of variable, we obtain 
\[
Q(\kappa)=\kappa\Big[ E\Big( \frac{1}{\kappa}\big) - \big( 1-\frac{1}{\kappa^2} \Big)K\Big( \frac{1}{\kappa} \Big) \Big]=\kappa H\Big(\frac{1}{\kappa}\Big)\,, 
\]
where $K$ and $E$ are complete elliptic integrals of first and second kinds. Since $Q(\kappa)=\frac{\kappa}{s}$, by the reciprocal identify \cite[19.7.3]{DLMF}, we obtain  $H\big(\frac{1}{\kappa}\big)=\frac{1}{s}$.

With the above preparation, by a precise control of $\chi_k(c)$ \cite[Theorem 1]{bonami2016uniform}:
\[
(\pi NW)^2\frac{1}{\Phi\big(\frac{2NW}{k} \big)^2}< \chi_k(c) < (\pi NW)^2\frac{1}{\Phi\big( \frac{2NW}{k+1}\big)^2}\,,
\]
we have
\begin{align}\label{proof lower bound of c2-chik using precise formula}
(\pi NW)^2\Big( 1- \frac{1}{\Phi\big(\frac{2NW}{k+1} \big)^2}\Big) < c^2-\chi_k< (\pi NW)^2 \Big(1-\frac{1}{\Phi\big(\frac{2NW}{ k} \big)^2} \Big)\,.
\end{align}
Note that \cite[Proposition 2]{bonami2016uniform} is a direct consequence by a global bound of $\Phi$.

To obtain a sharper lower bound of $B_k-A$ when $\frac{3}{2}NW<k<2NW$, we approximate $1- \frac{1}{\Phi^2}$ when $s\in\big[1, \frac{4}{3}\big]$. Note that $1- \frac{1}{\Phi^2}$ has a nice form, and can be bounded by two quadratic polynomials $R_L$ and $R_U$. By imposing $R_L(1)= 0= R_U(1)$ and $R_L(\frac{4}{3})= 1- \frac{1}{\Phi^2(\frac{4}{3})}= R_U(\frac{4}{3})$ and controlling the maximal deviation, we have 
\[
R_L(s) = 0.694 s^2 - 1.08 s + 0.386\ \mbox{ and }\ R_U(s)=-0.155 s^2 + 0.901 s - 0.746
\] 
so that $R_L(s)\leq 1- \frac{1}{\Phi^2}(s)\leq R_U(s)$ for all $s\in\big[1, \frac{4}{3}\big]$ and
\[
\Big\|1- \frac{1}{\Phi^2}-R_U\Big\|_{L^\infty(0,\frac{4}{3})}\leq 0.004\ \mbox{ and }\ \Big\|1- \frac{1}{\Phi^2}-R_L\Big\|_{L^\infty(0,\frac{4}{3})}\leq 0.022\,.
\]
While this is not a sharp bound, it is sufficient for our purpose.
As a result, we have
\begin{align}\label{proof lower bound of c2-chik using approximation formula}
&(\pi NW)^2
\Big( 
0.694 \Big(\frac{2NW}{k+1} \Big)^2 - 1.08 \Big(\frac{2NW}{k+1} \Big) + 0.386
\Big) < c^2-\chi_k\\
< \,& (\pi NW)^2 \Big( 
-0.155 \Big(\frac{2NW}{k} \Big)^2 +0.901 \Big(\frac{2NW}{k} \Big) - 0.746
\Big)\,.\nonumber
\end{align}

By the same control \cite[Proposition 2]{bonami2016uniform}, we have $\frac{\pi}{2}ck\leq \chi_k\leq 2c(k+1)$ for all $k<2NW$. Therefore, when $N$ is sufficiently large, $1-B_k=\frac{2\chi_k(c)}{(N-1)^2}(1+O(W^2))$ is controlled by $\frac{\pi^2 kW}{N-1}\leq 1-B_k\leq \frac{4\pi (k+1)W}{N}$. When $\frac{3}{2}NW<k<2NW$, we follow the same idea to get the claim. Indeed, consider 
\[
\tilde R_L(s) = 0.155 s^2 - 0.901 s + 1.746\ \mbox{ and }\ \tilde R_U(s)= -0.694 s^2 + 1.08 s + 0.614
\] 
so that $\tilde R_L(s)\leq \frac{1}{\Phi^2}(s)\leq \tilde R_U(s)$ for all $s\in\big[1, \frac{4}{3}\big]$ with $\tilde R_L(1)= 1= \tilde R_U(1)$ and $\tilde R_L(\frac{4}{3})= \frac{1}{\Phi^2(\frac{4}{3})}= \tilde R_U(\frac{4}{3})$. We have 
\[
\Big\|\frac{1}{\Phi^2}-\tilde R_U\Big\|_{L^\infty(0,\frac{4}{3})}\leq 0.022\ \mbox{ and }\ \Big\|\frac{1}{\Phi^2}-\tilde R_L\Big\|_{L^\infty(0,\frac{4}{3})}\leq 0.004\,, 
\]
which is sufficient for our purpose.
 
\end{proof}

The following lemma regards the size of $U_k(0)$. Recall that by Theorem \ref{theorem:said2023non:summary}(c), $0$ is not the maximal point of $|U_k|$ when $k>0$. We start with a remark.

\begin{remark}
The control of $U_k(0)$ seems trivial due to the following two existing results. First, it is stated in \cite[(73)]{slepian1978pswfV} that when $W\to 0$, $N=\lfloor \frac{c}{\pi W}\rfloor$, and $c>0$, then for a given $k$, when $N\to \infty$,
\begin{align}
U_k(0)\sim \frac{1}{\sqrt{W}}\psi^c_k(0)\,,\label{slepian1978 point approximation no rate}
\end{align}
where $\sim$ indicates asymptotically the ratio of left-hand and right-hand sides converges to $1$ and $\psi^c_k$ is the $k$-th PSWF with the time-bandwidth parameter $c=\pi NW$.
Second, when $k\leq \lfloor (NW)^\alpha\rfloor$, where $\alpha<1/2$, we have \cite[(8.86)]{osipov2013prolate}
\[
\psi^c_k(0)= \phi^{\sqrt{c}}_k(0)+O(c^{4\alpha-2})=\frac{(-1)^{k/2}c^{1/4}\sqrt{k!}}{\pi^{1/4}2^{k/2}(k/2)!}+O(c^{4\alpha-2})\,,
\] 
where $\phi^{\sqrt{c}}_k(x)$ is the scaled Hermite function of order $0$.
When $\lfloor (NW)^\alpha\rfloor<k<2NW$, by \cite[Theorem 2.1]{bonami2014uniform}, we know $|\psi^c_k(0)|\leq 1$.
As a result, when $k\leq \lfloor (NW)^\alpha\rfloor$, we obtain
\begin{equation}
U_k(0)\sim  
\frac{N^{1/4}}{W^{1/4}}\frac{\sqrt{k!}}{2^{k/2}(k/2)!}+O(N^{4\alpha-2}W^{4\alpha-\frac{5}{2}})\,,\label{control of Uk0 using pointwise control existing results}
\end{equation}
and when $\lfloor (NW)^\alpha\rfloor<k<2NW$,
\[
U_k(0)\sim \frac{1}{\sqrt{W}}|\psi^c_k(0)|\leq W^{-1/2}\,.
\]
Note that when $k$ is large, by Stirling formula, we have the approximation that
$\frac{\sqrt{k!}}{2^{k/2}(k/2)!}=\frac{2^{1/4}}{(\pi k)^{1/4}}\big(1-\frac{1}{8k}+O(k^{-2})\big)$,
which gives the decay $k^{-1/4}$ when $k\leq \lfloor (NW)^\alpha\rfloor$ is large.
In this work, we need a uniform bound on $U_k(0)$ over $k$, together with a convergence rate, so the above argument is insufficient. To the best of our knowledge, the literature does not provide an approximation rate that is uniform in $k$. 
\end{remark}

\begin{lemma} \label{taper_upper bound in frequency domain at 0, post-plunge}
Assume $W\to 0$ and $NW\to\infty$ when $N\to \infty$. 
When $N$ is sufficiently large, for any even $k=2,4,\ldots,2\lfloor \frac{N}{2}\rfloor$, we have
\[
|U_k(0)| \le 2\Big(\frac{1}{\vartheta_k}\Big)^{1/4}\sqrt{ \frac{\lambda_k }{ W }}\,,
\]
where $\vartheta_k:=\frac{k}{2NW}\wedge 1$.
\end{lemma}

\begin{proof}
Recall that $w_k:=\frac{1}{2\pi}\cos^{-1}(B_k)$ is the turning point for the DPSWF ODE \eqref{2nd order ODE for DSPWF}.  
By Theorem \ref{theorem:said2023non:summary}(a), when $k<2NW$, $w_k<W$, and we consider two cases. When $cNW<k<2NW$ for a fixed small constant $c>0$, by Lemma \ref{taper_upper bound in frequency domain, control of 1-A/B-A} and \eqref{proof of Ukf decay, wk and W in Bk and A}, we know 
\[
\sqrt{\frac{\vartheta_k}{2}\frac{N}{N-1}}W(1+O(W^2))<w_k<\sqrt{\frac{2\vartheta_k}{\pi}}W(1+O(W^2))\,,
\]
and hence set $\varrho_k=\frac{1}{2}\sqrt{\frac{\vartheta_k}{2}}$. When $0\leq k\leq cNW$, set $\varrho_k=\frac{1}{4}\sqrt{\frac{\vartheta_k}{2}}$. 
When $k\geq 2NW$, $w_k>W$, and we set $\varrho_k=1/2$. 
Below, due to the symmetry, all analyses are done over $I_k:=[0,\varrho_k W]$. 

Consider the Liouville-Green transformation (or WKB approximation). Set 
\begin{equation}\label{control Psik(0) when k>2NW using WKB, conversion0}
\Psi_k(f)=\sqrt{p(f)}U_k(f)\,, 
\end{equation}
and obtain
\begin{equation}
\Psi''_k(f)+Q_k(f)\Psi_k(f)=0\,,\label{control Psik(0) when k>2NW using WKB, equation of Psik}
\end{equation}
where 
\begin{equation}
Q_k(f):=\underbrace{4\pi^2 M\frac{q_k(f)}{p(f)}}_{:=E_1(f)}+\underbrace{\pi^2\frac{\cos^2(2\pi f)-2A\cos(2\pi f)+1}{p(f)^2}}_{:=E_0(f)}\,.\label{proof WKB Qk equation}
\end{equation}
Note that $Q_k=
\frac{\pi^2}{p(f)^2}[4Mp(f)q_k(f)+p(f)^2+(1-A^2)]$ is positive and $\sqrt{Q_k}$ can be viewed as the {\em instantaneous frequency} of the oscillatory solution $\Psi_k$ on $I_k$. 
Using the fact that $\frac{2}{\pi}x\leq \sin x\leq x$ when $x\in [0,\pi/2]$, the denominator of $E_1$ satisfies 
\[
2(1-\varrho_k^2)W^2\leq \cos(2\pi \varrho_k W) - A \leq p(f) \leq 1-A\leq 2\pi^2W^2\,, 
\]
and the numerator satisfies
\[
6\varrho_k^2 W^2= 8(w_k^2-f^2)\leq  q_k(f) \leq 1-B_k\leq 2\,. 
\]
Therefore, for any $f\in[0,\varrho_k W]$ and $k< 2NW$, we have
\[
3\varrho_k^2 M\leq E_1(f)\leq \frac{\pi^2}{1-\varrho_k^2}MW^{-2}\,.
\]
For $E_0(f)$, we have $\cos^2(2\pi f)-2A\cos(2\pi f)+1\leq 2(1-A)+A(1-\cos(2\pi \varrho_k W))\leq (4+2\varrho_k^2)\pi^2 W^2$, and hence $|E_0(f)|=O(W^{-2})$. Therefore, 
\[
\frac{|E_0(f)|}{|E_1(f)|}=O((NW)^{-1})\,, 
\]
where the worst bound happens when $k=2$, or $\rho_2=\sqrt{\frac{1}{2NW}}$. Hence, $E_0$ is dominated by $E_1$ under the assumption that $NW\to \infty$, and 
\[
Q_k(f)=E_1(f)(1+O((NW)^{-1})
\] 
for all $k$. 
Also, $Q_k(0)=\Theta( \frac{N}{W})$ when $k$ is small and $Q_k(0)=\Theta(N^2)$ when $k>cNW$.

Next, set the {\em phase} as 
\[
\xi(f):=\int^f_0\sqrt{Q_k(t)}dt,
\]
or $d\xi=\sqrt{Q_k(f)}df$ and factorize $\Psi_k$ into 
\begin{equation}
\Psi_k(f)=Q_k(f)^{-1/4}Y(\xi(f))\,.\label{control Psik(0) when k>2NW using WKB, conversion}
\end{equation}
By plugging \eqref{control Psik(0) when k>2NW using WKB, conversion} into \eqref{control Psik(0) when k>2NW using WKB, equation of Psik}, as a function of $\xi$, $Y(\xi)$ satisfies
\begin{equation}\label{control Psik(0) when k>2NW using WKB, conversion2}
\frac{d^2Y}{d\xi^2}+(1-\epsilon_k(\xi))Y(\xi)=0\,,
\end{equation}
where 
\[
\epsilon_k(\xi):=\epsilon_k(\xi(f))=\frac{5Q_k'(f)^2-4Q_k(f)Q_k''(f)}{16Q_k^3(f)}\,.
\]
By a direct calculation, we have
\begin{align*}
5Q_k'(f)^2-4Q_k(f)Q_k''(f)
=-\frac{4\pi^6[M^2B_2(f)+MB_1(f)+B_0(f)]}{p(f)^6}\,,
\end{align*}
where  
\begin{align*} 
B_0(f) =  4 (1 - A^2) \big[ -4 p(f)^4 - 10 A p(f)^3 + 7 (1 - A^2) p(f)^2 + (1 - A^2)^2 \big] \,,
\end{align*}
\begin{align*} 
B_1(f) =\,& 16 p(f)\big[ -(A - B_k) p(f)^4 - \big( 3 A (A - B_k) + 4 (1 - A^2) \big) p(f)^3 \\
 & - 2 (1 - A^2) (4 A + B_k) p(f)^2 + 3 (1 - A^2) \big( 2 (1 - A^2) - A (A - B_k) \big) p(f) \\ 
 & + 3 (A - B_k)(1 - A^2)^2  \big] \,,
\end{align*}
and
\begin{align*}
B_2(f) = \,&-16 (A - B_k) p(f)^2\big[ 4 p(f)^3 + (11 A + B_k) p(f)^2 \\
 & + 2 \big( A(A - B_k) - 4(1 - A^2) \big) p(f)- 3 (A - B_k)(1 - A^2)  \big] \,.
\end{align*}
As a result, 
\[
|\epsilon_k(\xi)|\leq \frac{|M^2B_2(f)+MB_1(f)+B_0(f)|}{256M^3p(f)^3q(f)^3(1+O((NW)^{-1})} \,.
\]
We comment that by a direct bound, $|\epsilon(\xi)|=O((NW)^{-2}\varrho_k^{-6})$ when $\varrho_k\to 0$ when $N\to \infty$. In other words, when $k$ is small, $\epsilon(\xi)$ can be as large as order $NW$.

As a perturbation of the equation $\frac{d^2Y}{d\xi^2}+Y(\xi)=0$, the exact solution to the transformed ODE \eqref{control Psik(0) when k>2NW using WKB, conversion2}, denoted as $Y_k$, satisfies the Volterra integral equation
\begin{equation}
    Y_k(\xi) = A_0 \cos(\xi + \phi) + \int_0^\xi \sin(\xi-v) \epsilon_k(v) Y_k(v) dv\,,
\end{equation}
where $A_0>0$ and $\phi\in [0,2\pi)$ are constants depending on the initial conditions of the transformed ODE at $0$. 
Clearly, $|Y_k(\xi)|\leq |A_0|+\int_0^{\xi}|\epsilon_k(v)||Y_k(v)|dv$. By Gronwall's inequality, we get the control of amplitude of $Y$:
\[
|Y_k(\zeta)|\leq |A_0|\exp\big(\int_0^{\zeta}|\epsilon_k(u)|du\big)=:|A_0|e^{V_k}\,.
\]
where $\zeta\in[0,\xi(\varrho_k W)]$ and
\[
V_k :=  \int_0^{\zeta}|\epsilon_k(u)|du=\int_0^{\varrho_k W}|\epsilon_k(\xi(f))|\sqrt{Q_k(f)}df
\]
by a change of variable $u=\xi(f)$.
By a direct bound with plugging in $\epsilon_k$ and $Q_k$, we have 
\begin{equation*}
    V_k  \leq \frac{20(1-\varrho_k^2)}{\sqrt{3}\pi \varrho_k^2(4-\varrho_k^2)^{1.5}NW}+O((NW)^{-3})\leq \frac{20(1-\varrho_k^2)}{\pi \varrho_k^2(4-\varrho_k^2)^{1.5}NW} 
\end{equation*}
when $\varrho_k$ is bounded away from $0$,
and
\begin{equation*}
    V_k  \leq  \frac{19\sqrt{15}(1-\varrho_k^2)}{540 \pi \varrho_k^2 NW}+O\Big(\frac{1}{NW}\Big)<0.044 \frac{1-\varrho_k^2}{\varrho_k^2 NW} 
\end{equation*}
when $\varrho_k\to 0$ as $N\to \infty$.
Note that the worst bound happens when $k=2$, where $V\leq 0.088$.
Therefore, for $\zeta\in[0,\xi(\varrho_k W)]$,  
\begin{align}
|Y_k(\zeta)-A_0 \cos(\zeta + \phi)|\,& \leq \int_0^{\zeta}|\sin(\zeta-u)||\epsilon_k(u)||Y_k(u)|du\nonumber\\
&\leq  |A_0|e^{V_k}\int_0^{\zeta} |\epsilon_k(u)|du = |A_0|V_ke^{V_k}\,.\label{proof of V_k(0) upper bound for postplunge, Y control-pre}
\end{align}
Since the worst case happens as $V_2e^{V_2} \leq 0.0961$, for any even $k$ we have
\begin{align}
Y_k(\zeta)=A_0(\cos(\zeta+\phi)+e_k(\zeta))\,,\label{proof of V_k(0) upper bound for postplunge, Y control}
\end{align}
where $\|e_k\|_{L^\infty[0,\xi(\varrho_k W)]}\leq \|e_2\|_{L^\infty[0,\xi(\varrho_k W)]}\leq 0.0961$.

Next, recall the relationships \eqref{control Psik(0) when k>2NW using WKB, conversion0} and \eqref{control Psik(0) when k>2NW using WKB, conversion} that $U_k^2(f) = E_k(f) |Y_k(\xi(f))|^2$, where the deterministic envelope is 
\[
E_k(f) := [Q_k(f)^{1/2} p(f)]^{-1}\,.
\]
Expanding the dominant term of $Q_k$ yields 
\[
Q_k(f)^{1/2} p(f) =\pi N \sqrt{p(f)q_k(f)} +O(1)\,,
\]
which is strictly decreasing on $[0, \varrho_k W]$. Hence, $E_k(f)$ is strictly and monotonically increases on $[0, \varrho_k W]$, guaranteeing $E_k(f) \geq E_k(0)$ for any $f\in [0,\varrho_k W]$.

We now put all the above together. By a direct bound, we have
\begin{align*}
    \lambda_k \,&=  \int_{- W}^{ W} U_k^2(f) df > 2 \int_0^{\varrho_k W} E_k(f) |Y_k(\xi(f))|^2 df
    \\
    &> 2 E_k(0) A_0^2 \int_0^{\varrho_k W} (\cos^2(\xi(f) + \phi)-\|e_k\|_\infty |\cos(\xi(f)+\phi)|) df \,,
\end{align*}
where we use the monotonicity $E(f) > E(0)$ and a direct lower bound of \eqref{proof of V_k(0) upper bound for postplunge, Y control} as 
\begin{align*}
|Y_k(\zeta)|^2\,&=|A_0(\cos(\zeta+\phi)+e_k(\zeta))|^2 \geq A_0^2(\cos(\zeta+\phi)^2-\|e_k\|_\infty)_+\\
&\geq A_0^2\cos(\zeta+\phi)^2-2\|e_k\|_\infty A_0|\cos(\zeta+\phi)|\,.
\end{align*}
%
By trigonometric identity, $\int_0^{\varrho_k W} \cos^2(\xi(f) + \phi) df$ becomes
\begin{equation*}
    \int_0^{\varrho_k W} \frac{1 + \cos(2\xi(f)+2\phi)}{2} df = \frac{\varrho_k W}{2} + \frac{1}{2} \int_0^{\varrho_k W} \cos(2\xi(f)+2\phi) df\,.
\end{equation*}
To control the size of $\int_0^{\varrho_k W} \cos(2\xi(f)+2\phi) df$, since $\frac{d}{df}[2\xi(f)+2\phi]=2\sqrt{Q_k(f)}$, by integration by parts  we have
\begin{align*}
    \int_0^{\varrho_k W}&\cos(2\xi(f)+2\phi) df =  \int_0^{\varrho_k W} \frac{1}{2\sqrt{Q_k(f)}}\frac{d}{df}[\sin(2\xi(f)+2\phi)] df  \\
    &=\frac{\sin(2\xi(f)+2\phi)}{2\sqrt{Q_k(f)}} \Big|_0^{\varrho_k W} - \int_0^{\varrho_k W} \sin(2\xi+2\phi) \frac{d}{df} \Big(\frac{1}{2\sqrt{Q_k(f)}}\Big) df\,.
\end{align*}
We have
\[
\left|\frac{\sin(2\xi+2\phi)}{2\sqrt{Q_k(f)}} \Big|_0^{\varrho_k W}\right|\leq \frac{1}{2\sqrt{Q_k(\varrho_k W)}} +\frac{1}{2\sqrt{Q_k(0)}}\,. 
\]
Because $Q_k(f)$ is strictly increasing, $1/\sqrt{Q_k(f)}$ is strictly decreasing, and hence
\begin{align*}
&\Big|\int_0^{\varrho_k W} \sin(2\xi+2\phi) \frac{d}{df} \Big(\frac{1}{2\sqrt{Q_k(f)}}\Big) df\Big|\leq \int_0^{\varrho_k W} \Big| \frac{d}{df} \Big(\frac{1}{2\sqrt{Q_k(f)}}\Big)\Big| df\\
=\,&-\int_0^{\varrho_k W}  \frac{d}{df} \Big(\frac{1}{2\sqrt{Q_k(f)}}\Big) df=-\frac{1}{2\sqrt{Q_k(\varrho_k W)}} +\frac{1}{2\sqrt{Q_k(0)}}\,.
\end{align*}
As a result,
\[
\Big|\int_0^{\varrho_k W} \cos(2\xi(f)+2\phi) df\Big|\leq \frac{1}{\sqrt{Q_k(0)}}\,. 
\]
Thus, 
\[
\int_0^{\varrho_k W} \cos^2(\xi(f) + \phi) df\leq \frac{\varrho_k W}{2} +\frac{1}{\sqrt{Q_k(0)}}\,.
\] 
Apply a trivial upper bound, we have $\int_0^{\varrho_k W} \|e\|_\infty |\cos(\xi(f)+\phi)| df\leq \varrho_k W \|e\|_\infty$.
By Lemma \ref{taper_upper bound in frequency domain, control of 1-A/B-A}, 
\[
Q_k(0)=4\pi^2 M\frac{1-B_k}{1-A}\Big(1+O\Big(\frac{1}{NW}\Big)\Big)\leq \frac{2\pi  (k+1)N}{W}(1+O(W^2))\,, 
\]
which leads to 
\[
\frac{1}{\sqrt{Q_k(0)}}\geq \frac{k+1}{\pi k}\varrho_k W\,.
\] 
As a result, we have 
\begin{equation*}
    \lambda_k > 2E(0) A_0^2 \Big( \varrho_k W\Big(\frac{1}{2}-\|e_k\|_\infty\Big) +\frac{1}{\sqrt{Q_k(0)}} \Big)> 2E(0) A_0^2  \varrho_k W\Big(\frac{1}{2}-\|e_k\|_\infty+\frac{k+1}{\pi k}\Big) \,.
\end{equation*}
On the other hand, at $0$ we have
\[
U_k^2(0) = E(0) |Y(0)|^2 \leq E(0) A_0^2 (1 + \|e_k\|_\infty)\,.
\]
By the above two bound, we get 
\begin{equation*}
    |U_k(0)| \leq \sqrt{ \frac{\lambda_k(1+\|e_k\|_\infty)}{\varrho_k W(\frac{1}{2}-\|e_k\|_\infty+\frac{k+1}{\pi k})}} \,.
\end{equation*}
The claim comes from plugging back $\varrho_k$ and use the worst control; that is, when $k=2$, we have $\sqrt{2\sqrt{2}\frac{1+\|e_2\|_\infty}{\frac{1}{2}-\|e_2\|_\infty+\frac{3}{2\pi}}}\leq 2$.
\end{proof}

\begin{lemma} \label{taper_upper bound in frequency domain, pre-plunge}
Assume $W\to 0$ and $NW\to\infty$ when $N\to \infty$. 
When $0\leq k<  2NW$ is even, we have 
\[
 |U_k(1/2)| \leq 2\sqrt{1-\lambda_k} (1+O(N^{-1}))
  \]
 when $N$ is odd, and $|U_k(1/2)|=0$ when $N$ is even.
\end{lemma}

\begin{proof}

The proof is almost the same as that for Lemma \ref{taper_upper bound in frequency domain at 0, post-plunge} but much simpler since $p$ and $q$ have the same sign and do not degenerate over $[1/4,1/2]$ under assumption. We sketch the key steps here. Set $\Psi_k(f)=\sqrt{p(f)}U_k(f)$ over $[W,1/2]$ and $\xi(f):=\int^f_W\sqrt{Q_k(t)}dt$, where $Q_k$ is defined in \eqref{proof WKB Qk equation}, the Liouville-Green transform gives 
\begin{equation}\label{control Uk(1/2) when k>2NW using WKB, conversion2}
\frac{d^2Y}{d\xi^2}+(1-\epsilon(\xi))Y(\xi)=0\,,
\end{equation}
where $|\epsilon(\xi)|=O(N^{-2})$. Therefore, 
\begin{align*}
|Y(\zeta)-A_0 \cos(\zeta + \phi)| =|A_0|O(N^{-1})\,,
\end{align*}
where the implied constant is independent of $k$, and $A_0$ and $\phi$ depends on the initial conditions of the transformed ODE at $1/2$. Since
\begin{align*}
    1-\lambda_k =\,&  \int_{[-\frac{1}{2},\frac{1}{2}]\backslash (-W,W)} U_k^2(f) df > 2 \int_{1/4}^{1/2} E(f) |Y(\xi(f))|^2 df\\
    >\,& 2 E(1/2) A_0^2 \int_{1/4}^{1/2} [\cos^2(\xi(f) + \phi)+O(N^{-1})] df\\
    =\,& 2 E(1/2) A_0^2 \Big( \frac{1}{8}+\frac{1}{\sqrt{Q_k(1/2)}}+O(N^{-1}) \Big)=E(1/2) A_0^2 \Big( \frac{1}{4}+O(N^{-1}) \Big)\,,
\end{align*}
where $E(f) := [Q_k(f)^{1/2} p(f)]^{-1}$ is the envelope and $Q_k(1/2)=\Theta(N^2)$. Note that $Q_k(f)>0$ here is strictly increasing from $1/4$ to $1/2$, while $E(f)$ is strictly decreasing.
At $1/2$, we have
\[
U_k^2(1/2) = E(1/2) |Y(1/2)|^2 \leq E(1/2) A_0^2 (1 + O(N^{-1}))\,.
\]
Jointly, when $N$ is sufficiently large, for any even $k< 2NW$ we have the claim
\begin{equation*}
    |U_k(1/2)| \leq  2 \sqrt{ 1-\lambda_k } (1+O(N^{-1}))\,.
\end{equation*}

\end{proof}

\begin{lemma} \label{taper_upper bound in frequency domain, derivative at 0 and 1/2 control}
Assume $W\to 0$ and $NW\to\infty$ when $N\to \infty$. 
When $k< 2NW$ is even and $N$ is odd, $U'_k(1/2)=0$. When $k< 2NW$ is even and $N$ is even, we have
\begin{equation}
    |U_k'(1/2)|^2 \leq 7\pi^2 (N^2-1)(1-\lambda_k)\,. 
\end{equation}
\end{lemma}

The $N^2-1$ factor is natural for derivative control. When we apply this bound, the $N^2-1$ factor is canceled by normalization; for example, \eqref{application of Said 2023 control for outside band for k<2NW}.

\begin{proof}
Motivated by the proof for \cite[Theorem 2]{said2023non}, construct an auxiliary functional $Q(f)$ by
\begin{align}
    Q(f) &\,= U_{k}^2(f) + \frac{p(f)}{4\pi^2 Mq_k(f)} U_{k}'(f)^2\label{proof U'k(0) pretplunge, Q function control}\,.
\end{align}
Similar to the proof of \cite[Theorem 2]{said2023non}, by differentiating this with respect to $f$ and substituting $U_k''(f)$ 
using \eqref{2nd order ODE for DSPWF} yields cross-term cancellation and gives \cite[(51)]{said2023non} we have
\begin{equation*}
    Q'(f) = \frac{ \sin(2\pi f) (p(f)+q_k(f))}{2\pi Mq_k(f)^2} U_{k}'(f)^2\,,
\end{equation*}
which is strictly negative over $(W,1/2]$ except finite points (see Theorem \ref{theorem:said2023non:summary}(f)) since when $k < \lfloor 2NW\rfloor$, we have $B_k > A$. Therefore, $Q(f)$ is monotonically decreasing on $(W,1/2]$ and we have a trivial bound 
\[
Q(1/2)\leq \frac{1}{\epsilon} \int_{1/2-\epsilon}^{1/2} Q(f)df
\] 
for any $\epsilon\in(0,1/2-W)$. 

Next, we claim that when $N$ is sufficiently large, we can find $t_k\in [1/5,1/4]$ that is a root of $U_k$. 
Note that in the Liouville-Green transform \eqref{control Psik(0) when k>2NW using WKB, equation of Psik}, we have $\Psi''_k(f)+Q_k(f)\Psi_k(f)=0$, where $\Psi_k(f)=\sqrt{p(f)}U_k(f)$ and $Q_k$ is defined in \eqref{proof WKB Qk equation}, and $Q_k=\Theta(N^2)$ over $[1/5,1/2]$.
By the root distance control by \cite[Theorem 2.18]{osipov2013prolate}, the distance of two consecutive roots of $\Psi_k(f)$ is of order $\Theta(N^{-2})$. Since $p(f)$ is positive over $[1/5,1/2]$, all roots of $\Psi_k$ are roots of $U_k$, and we obtain the claim.

When $k$ and $N$ are both even, $U_k(1/2)=0$. Taking $t_k<1/4$ that is the closest root of $U_k$ to $1/4$. Plug $f=1/2$ into \eqref{proof U'k(0) pretplunge, Q function control} and obtain
\begin{align}
&\frac{1+ A}{4\pi^2 M(1+B_k)} U_{k}'(1/2)^2\leq
4\int_{t_k}^{1/2}U_{k}^2(f) + \frac{p(f)}{4\pi^2 Mq_k(f)} U_{k}'(f)^2df \label{control of U_k' over t_k and 1/2, step0} \\
\leq \,&2(1-\lambda_k)+ \int_{t_k}^{1/2}\frac{p(f)}{\pi^2 Mq_k(f)} U_{k}'(f)^2 d\leq 2(1-\lambda_k)+ \frac{1}{\pi^2M}\int_{t_k}^{1/2} U_{k}'(f)^2 df\,,\nonumber
\end{align}
where we use $\frac{p(f)}{q_k(f)}\leq 1$ over $f\in[t_k,1/2]$ when $N$ is sufficiently large. Since $U_k(t_k)=U_k(1/2)=0$, by multiplying $U_k(f)$ to the ODE of DPSWF, integrating, and applying integration by part, we have
\[
\int_{t_k}^{1/2}(-p(f))U_k'(f)^2df=4\pi^2 M \int_{t_k}^{1/2}(-q_k(f))U_k^2(f)df\,.
\]
Since $-p(f)>1$ and $-q_k(f)\leq 2$ on $[t_k,1/2]$, we have 
\begin{align}
\int_{t_k}^{1/2}U_k'(f)^2df\leq 8\pi^2M \int_{t_k}^{1/2}U_k^2(f)df\leq 4\pi^2M(1-\lambda_k)\label{control of U_k' over t_k and 1/2}
\end{align}
when $N$ is sufficiently large since $-p$ and $-q_k$ are both strictly positive. By plugging \eqref{control of U_k' over t_k and 1/2} to \eqref{control of U_k' over t_k and 1/2, step0}, we have
\begin{align*}
 U_{k}'(1/2)^2\leq 7\pi^2 (N^2-1)(1-\lambda_k)\,,
\end{align*}
where we use $M=(N^2-1)/4$ and the bound that $\frac{1+B_k}{1+ A}\leq \frac{7}{6}$ when $N$ is sufficiently large, and hence the claim.
\end{proof}

\begin{lemma} \label{taper_upper bound in frequency domain, f>W k<2NW case}
Assume $W\to 0$ and $NW\to\infty$ when $N\to \infty$. When $N$ is sufficiently large,  for any even $k< 2NW$, we have
\begin{align}
|U_k(f)|\leq  \frac{4\sqrt{1 - \lambda_k}}{\sqrt{p(f)q_k(f)}}\,,\label{bound Uk over W,1/2 for preplung}
\end{align}
for $f\in (W,1/2]$. Moreover,
\begin{equation}
    |U_k(f)| \leq  2N^2\sqrt{1 - \lambda_k}\Big(1+O(W^2)+O\Big(\frac{1}{N^2W^2}\Big)\Big)\label{bound Uk over W,1/2 for preplung2}
\end{equation}    
for $f\in [W, (1+\nu_k)W]$, where $\nu_k$ satisfies $\nu_k=\frac{\sin(2\pi W)}{\pi W(N^2-1)(B_k-A)}=\Theta\big(\frac{1}{N^2W^2}\big)$. 
\end{lemma}

When $f>(1+\tau)W$ for any $\tau>0$, the denominator of \eqref{bound Uk over W,1/2 for preplung} is bounded from below, so that $|U_k(f)|$ is of order $\sqrt{1-\lambda_k}$. The control over $f\in [W, (1+\nu_k)W]$ is worst when $k\sim 2NW$, since $\lambda_k$ is no longer exponentially close to $1$. In this case, $B_k-A$ is small, and $\nu_k$ is large. 

\begin{proof}
Consider the envelope function 
\[
E(f)=\pi^2(N-1)^2 p(f)q_k(f)U_k(f)^2+(p(f)U_k'(f))^2\geq 0\,. 
\]
Below we drop $f$ to ease the notation. By a direct calculation, we have $(pU_k')'=-q_kU_k$, and hence $E'=\pi^2(N-1)^2(pq)'U_k^2\geq 0$ since $p\leq 0$ and $q<0$ are both monotonically decreasing over $(W,1/2]$. Therefore, $E$ is increasing. 
When $N$ is odd, by Lemma \ref{taper_upper bound in frequency domain, pre-plunge} we have 
\begin{align*}
E(f)\leq \pi^2(N-1)^2p\Big(\frac{1}{2}\Big)q_k\Big(\frac{1}{2}\Big)U_k\Big(\frac{1}{2}\Big)^2 \leq 16\pi^2 (N-1)^2(1 - \lambda_k)\,.
\end{align*}
where we use $p\big(\frac{1}{2}\big)q_k\big(\frac{1}{2}\big)<4$.
When $N$ is even, by Lemma \ref{taper_upper bound in frequency domain, derivative at 0 and 1/2 control} we have 
\[
E(f)\leq p\Big(\frac{1}{2}\Big)^2U'_k\Big(\frac{1}{2}\Big)^2\leq 14\pi^2 (N-1)^2 (1 - \lambda_k)\,.
\]
Since $\pi^2(N-1)^2 pqU_k^2\leq E$, we conclude the first claim \eqref{bound Uk over W,1/2 for preplung} over $(W,1/2]$ that
\[
|U_k(f)|\leq  \frac{4\sqrt{1 - \lambda_k}}{\sqrt{p(f)q_k(f)}}\,.
\]

The bound of \eqref{bound Uk over W,1/2 for preplung} blows up at $W$, so we control it differently. Note that $p$ is degenerate at $W$ with a simple root, so the standard Frobenius asymptotic expansion shows that near $W$ we have
\[
U_k(f)=c_1u_1(f)+c_2u_2(f)\,,
\]
where 
\begin{align*}
u_1(f)=1+\sum_{l=1}^\infty a_l(f-W)^l
\end{align*}
is the regular solution with the radius of convergence $2W$ 
with
\[
a_1=-\frac{\pi (N^2-1)q_k(W)}{p'(W)}(f-W)  
\]
 and $u_2$ is the singular solution that blows us at $W$ with the rate $\ln(f-W)$ when $f\to W+$. 
Indeed, by plugging Taylor expansion of $p(f)=\sum_{l=1}^\infty p_l(f-W)^l$ and $q_k(f)=\sum_{l=0}^\infty q_{k,l}(f-W)^l$ near $W$ into the DPSWF's ODE, and matching the order of $(f-W)^l$, 
we get the indicial equation at $f=W$ as $r^2=0$, and hence the fundamental solutions, where $a_l$ satisfies the the recurrence
\[
a_{l+1}
=
-\frac{
(l+1)\sum_{j=2}^{l+1}
(l+2-j)p_j a_{l+2-j}
+\pi^2(N^2-1)\sum_{j=0}^{l}q_{k,j} a_{l-j}
}{
(l+1)^2p_1
}
\]
for any $l\geq 1$ and $a_0=1$.
Since $U_k$ is bounded at $W$ by Theorem \ref{theorem:said2023non:summary}(e), we conclude that $c_2=0$ and hence near $W$,
\[
U_k(f)=c_1u_1(f)= c_1\Big( 1+\sum_{l=1}^\infty a_l(f-W)^l \Big)\,.
\] 
Note that $u_1(W)=1$.

Next, we control $c_1$. Consider $f_{\nu_k}:=(1+\nu_k)W>W$, where $\nu_k>0$ is small depending on $N$ that we will determine later. 
By the first claim \eqref{bound Uk over W,1/2 for preplung}, we have
\[
U_k(f_{\nu_k})\leq  \frac{4\sqrt{1 - \lambda_k}}{\sqrt{p(f_{\nu_k})q_k(f_{\nu_k})}}\leq  \frac{2\sqrt{1 - \lambda_k}}{\pi^2\nu_k(2+\nu_k) W^2}(1+O(W^2))\,,
\]
where we use $p(f_{\nu_k})=-2\sin(\pi(2W+\nu_k W))\sin(\pi\nu_k W)$, $q_k(f_{\nu_k})=-2\sin(\pi((1+\nu_k)W+w_k))\sin(\pi((1+\nu_k)W-w_k))$, and hence $p(f_{\nu_k})q_k(f_{\nu_k})\geq 4\pi^4\nu_k^2(2+\nu_k)^2W^4+O(W^6)$. 
On the other hand, by integrating $(pu_1')=-\pi^2(N^2-1)qu_1$ twice with $u_1(W)=1$, we have
\[
u_1(f)=1-\pi^2(N^2-1)\int_W^f \frac{1}{-p(t)}\int_W^t (-q_k(s))u_1(s)dsdt\,.
\]
Since $u_1(W)=1$ and $u_1$ is analytic near $W$, $u_1$ is positive near $W$. As long as $u_1$ is positive, $\frac{1}{-p(t)}\int_W^t (-q_k(s))u_1(s)ds$ is positive, and hence $u_1$ is strictly monotonically decreasing after $W$. So there exists $f_*>W$ so that $u_1(f)\geq1/2$ for all $f\in [W, f_*]$. Over $[W,f_*]$, we have
\[
1-I(f)\leq u_1(f)\leq 1-\frac{1}{2}I(f)\,,
\]
where $I(f):=\pi^2(N^2-1)\int_W^f \frac{1}{-p(t)}\int_W^t (-q_k(s))dsdt$. Next, determine $\nu_k$ so that $I(f_{\nu_k})\geq  1/2$. Plugging Taylor expansion of $p$ and $q_k$ near $W$ into $I$, we obtain
\[
I(f)=|a_1|(f-W)\Big(1+O\Big(\frac{f-W}{W}\Big)\Big)\,,
\] 
and hence the problem is reduced to
\[
(f-W)\Big(1+O\Big(\frac{f-W}{W}\Big)\Big)=\frac{I(f)}{|a_1|}\leq \frac{1}{2|a_1|}=\frac{\sin(2\pi W)}{\pi (N^2-1)(B_k-A)}\,.  
\]
Choose $\nu_k=\frac{\sin(2\pi W)}{\pi W(N^2-1)(B_k-A)}$. By Lemma \ref{taper_upper bound in frequency domain, control of 1-A/B-A} and the fact that $B_0-A=2\pi^2W^2(1+O(W^2))$ since $B_2-A<B_0-A<1-A$, when $0\leq k\leq \frac{3}{2}NW$, we have
\[
\frac{1}{\pi^2N^2W^2}(1+O(W^2)) \leq \nu_k\leq  \frac{1}{(\pi^2-3\pi)N^2W^2}(1+O(W^2))\,,
\]
and when $\frac{3}{w}NW<k<2NW-1$, 
\[
\frac{1}{1.182N^2W^2}(1+O(W^2)) \leq \nu_k\leq \frac{1}{0.462N^2W^2}(1+O(W^2))\,. 
\]
Therefore, $O\big(\frac{f-W}{W}\big)=O\big(\frac{1}{N^2W^2}\big)$, and over $[W,(1+\nu_k)W]$, 
\[
I(f)= \frac{\pi(N^2-1)(B_k-A)}{2\sin(2\pi W)}\nu_kW\Big(1+O\Big(\frac{1}{N^2W^2}\Big)\Big)=\frac{1}{2}\Big(1+O\Big(\frac{1}{N^2W^2}\Big)\Big)\,.
\]
Putting the above together, we have
\[
\frac{|U_k(f_{\nu_k})|}{|u_1(f_{\nu_k})|}\leq \frac{4\sqrt{1 - \lambda_k}}{\pi^2\nu_k(2+\nu_k)W^2}\Big(1+O(W^2)+O\Big(\frac{1}{N^2W^2}\Big)\Big)\,.
\]
Since $c_1=\frac{U_k(f_{\nu_k})}{u_1(f_{\nu_k})}$, $u_1(W)=1$, and $\pi^2\nu_k(2+\nu_k)W^2>2\pi^2\nu_kW^2\geq \frac{2}{N^2}(1+O(W^2))$ for all $k<2NW$, the above controls give us  the bound
\[
|U_k(W)|\leq 2N^2\sqrt{1 - \lambda_k}\Big(1+O(W^2)+O\Big(\frac{1}{N^2W^2}\Big)\Big)\,.
\]
By the maximal principle, we obtain 
\[
|U_k(f)|\leq 2N^2\sqrt{1 - \lambda_k}\Big(1+O(W^2)+O\Big(\frac{1}{N^2W^2}\Big)\Big)
\]
for any $f\in [W, (1+\nu_k)W]$.

\end{proof}

\begin{lemma} \label{taper_upper bound in frequency domain, f<W k>2NW case}
Assume $W\to 0$ and $NW\to\infty$ when $N\to \infty$. Set $G:=\lfloor \frac{6}{\pi^2}\log(100NW+25)\log(N)\rfloor$ and fix $\upsilon<1$. When $ 2NW\leq k\leq 2NW+G$ is even and $N$ is sufficiently large, over $f\in [0,\upsilon W]$ we have
\[
|U_k(f)|\leq \frac{1}{1-\upsilon^2}\sqrt{\frac{\lambda_k}{W}}\,.
\]
\end{lemma}

Note that the same proof for Lemma \ref{taper_upper bound in frequency domain, f>W k<2NW case} can be extended to control $U_k(f)$ over the whole period $[0,W]$ and $k\geq 2NW$, as long as we generalize Lemma \ref{taper_upper bound in frequency domain, control of 1-A/B-A} to control $B_k-A$. However, this is not needed in our main proof, so we skip it.

\begin{proof}

The proof is almost the same as that of Lemma \ref{taper_upper bound in frequency domain, f>W k<2NW case}. The only difference is that $E'\leq 0$ since $p\geq 0$ and $q_k>0$. Therefore, we have 
\[
p(f)q_k(f)U_k(f)^2 \leq E(f)\leq E(0)=p(0)q_k(0)U_k(0)^2=4\sin(\pi W)^2 \sin(\pi w_k)^2U_k(0)^2\,. 
\]
By Lemma \ref{taper_upper bound in frequency domain at 0, post-plunge},  we obtain
\[
|U_k(f)|\leq \frac{4\pi w_k}{\sqrt{p(f)q_k(f)}}\sqrt{\lambda_kW}\,.
\]
Over $[0,\upsilon W]$ for $\upsilon<1$, we have $p(f)=2\pi^2(W^2-f^2)(1+O(W^2))\geq 2\pi^2W^2(1-\upsilon^2)(1+O(W^2))$. By Theorem \ref{theorem:said2023non:summary}(f), $w_k>W$ is an increasing sequence and $B_k>A$ is decreasing. Consider the worst bound $q_k(f)\geq p(f)\geq 2\pi^2W^2(1-\upsilon^2)(1+O(W^2))$, which leads to
\[
|U_k(f)|\leq \frac{2 w_k\sqrt{\lambda_kW}}{\pi W^2(1-\upsilon^2)}(1+O(W^2))\,.
\]
What remains is controlling $w_k$.
When $2NW\leq k\leq 2NW+G$, note that $\frac{1-B}{2\pi}=\frac{\sin(\pi w_k)^2}{\pi}\geq \frac{4w_k^2}{\pi}$. Denote $c:=\pi NW$. Note that $3<G<\frac{5\pi}{4}c$, which satisfies \cite[(3.21)]{osipov2013prolate}. So we can apply Theorem \ref{lemma: control M(1-Bk) asymptotically} and \cite[Theorem 3.7]{osipov2013prolate} to obtain 
\[
w_k^2\leq \frac{1-B_k}{8}\leq \frac{\chi_k(c)}{16M}\leq \frac{\pi^2W^2}{4} \Big(1+\frac{8G}{\pi^2 NW}\Big)
\]
when $N$ is sufficiently large.
Therefore, over  $[0,\upsilon W]$, we have
\[
|U_k(f)|\leq \frac{1}{1-\upsilon^2}\sqrt{\frac{\lambda_k}{W}}\,.
\]

\end{proof}

\section{Proof of Theorem \ref{theorem: bound Vf energy leakage uniform version}} \label{section proof of theorem theorem: bound Vf energy leakage uniform version}

We start with a well known fact about DPSWF with proofs for the sake of self-containedness.

\begin{lemma}
Suppose $W\in (0,1/2)$ and $N>1$. 
For $f\in [0,1/2]$, we have
\begin{equation}
\sum_{j=0}^{N-1}\lambda_j U_j(0)U_j(f)=\int_{-W}^W D_N(f-x)D_N(x)dx\,,\label{sum lamda VVf = conv of DN}
\end{equation}
where $D_N(x)=\frac{\sin(\pi Nx)}{\sin(\pi x)}$ is the Dirichlet kernel. We also have
\begin{align}
\int_{-1/2}^{1/2} D_N(x)^2dx=N\,.\label{int DN2 on I is N}
\end{align}
Third, for all $k=0,\ldots,N-1$, we have
\[
\sum_{j=0}^{N-1}U_j(0)\vv_{j,k}=1
\ \mbox{ and }\ 
\sum_{j=0}^{N-1}\lambda_j U_j(0)^2=\boldsymbol{1}^\top A\boldsymbol{1}\,,
\]
\end{lemma}
\begin{proof}
By \eqref{band_eigen}, we have $\int_{-W}^WD_N(f-x)U_j(x)dx=\lambda_j U_j(f)$. Therefore, 
\[
\sum_{j=0}^{N-1}\lambda_j U_j(0)U_j(f)=\int_{-W}^WD_N(f-x)\sum_{j=0}^{N-1}U_j(0)U_j(x)dx\,.
\]
By construction, we have 
\begin{align*}
\sum_{j=0}^{N-1} U_j(0)U_j(x)\,&=\sum_{m,n=0}^{N-1}e^{-i2\pi x(n-(N-1)/2)}\sum_{j=0}^{N-1} \vv_{j,m}\vv_{j,n}\\
&=\sum_{n=0}^{N-1}e^{-i2\pi x(n-(N-1)/2)}=D_N(x)\,,
\end{align*}
since $\sum_{j=0}^{N-1} \vv_{j,m}\vv_{j,n}=\delta_{m,n}$. The second claim comes from the Plancherel identity.
The third one comes from a direct algebraic evaluation. Since $\vv_k$ form an orthonormal basis of $\mathbb{R}^N$ and $U_k(0)=\boldsymbol{1}^\top\vv_k$, we have
\[
\boldsymbol{1}=\sum_{j=0}^{N-1}\boldsymbol{1}^\top \vv_k \vv_k=\sum_{j=0}^{N-1} U_k(0)\vv_k\,.
\]
The final equality comes from the eigenvalue decomposition of $A$ and $U_j(0)=\boldsymbol{1}^\top \vv_j$.
\end{proof}

The proof of Theorem \ref{theorem: bound Vf energy leakage uniform version} is split into several lemmas. By \eqref{sum lamda VVf = conv of DN} we have this ``oscillation cancellation'' identity:
\begin{align}
&\underbrace{D_N(f)}_{=:I_0}-\underbrace{\int_{|x|>W}D_N(f-x)D_N(x)dx}_{=:I_1}=\sum_{j=0}^KU_{j}(0)U_{j}(f)\nonumber\\
&\qquad\qquad+\underbrace{\sum_{j=0}^K(\lambda_{j}-1)U_{j}(0)U_{j}(f)}_{=:I_2}+\underbrace{\sum_{j=K+1}^{N-1}\lambda_{j}U_{j}(0)U_{j}(f)}_{=:I_3}\label{EQ:proof of lemma of sumUkUk(f) decomposition when f <W}
\end{align}
that converts the summation of oscillatory DPSWF $\sum_{j=0}^KU_{j}(0)U_{j}(f)$, the key quantity we have interest in the proof, into Dirichlet kernel operations, where we use the reproducing kernel property of $D_N$:
\[
\int_{-1/2}^{1/2}D_N(f-x)D_N(x)dx=D_N(f)=\frac{\sin((2N+1)\pi f)}{\sin(\pi f)}\,.
\]

\begin{lemma}\label{lemma summation of V_j(0)^2 is N}
Suppose $W=W(N)\to 0$ and $WN\to \infty$ as $N\to \infty$, and $K=\lfloor 2NW\rfloor-G$, where $G:=\lfloor \frac{2c}{\pi^2}\log(100NW+25)\log(N)\rfloor$ with $c>2$. We have
\[
\sum_{j=0}^K U_{j}(0)^2=N+O\Big(\frac{\log(NW)\log(N)}{W}\Big)
\]
when $N$ is sufficiently large.
\end{lemma}

\begin{proof}
We use \eqref{EQ:proof of lemma of sumUkUk(f) decomposition when f <W} with $f=0$. $I_0=N$ by a direct evaluation.
For $I_1$, since $f\in [W,1/2]$, we use the simple lower bound $2f\leq \sin(\pi f)\leq \pi f$ to simplify the integration into
\[
\int_{W}^{1/2}\frac{\sin(\pi Nf)^2}{\pi^2 f^2}df\leq \int_{W}^{1/2}D_N(f)^2df\leq \int_{W}^{1/2}\frac{\sin(\pi Nf)^2}{4f^2}df\,.
\]
By trigonometric identity and change of variable, we have 
\[
\int_{W}^{1/2}\frac{\sin(\pi Nf)^2}{f^2}df=\frac{\pi N}{2}\int_{\pi NW}^{\pi N/2}\frac{1-\cos(2x)}{x^2}dx=\frac{1}{2W}-\frac{1}{2}+O((N^2W)^{-1}),
\]
and hence $I_1=\frac{1}{2W}+O(1)$ since $NW\to \infty$. 
For $I_2$, we use \eqref{control of Uk0 using pointwise control existing results} to control $U_0(0)$ and Lemma \ref{taper_upper bound in frequency domain at 0, post-plunge} to control $U_j(0)$ when $2\leq j\leq K$ and obtain 
\[
I_2\leq (1-\lambda_0)\frac{N^{1/2}}{W^{1/2}}(1+O((NW)^{-9/4}))+\frac{4}{W}\sum_{j=2}^K\frac{1-\lambda_j}{\sqrt{\vartheta_j}} =O(N^{-1})\,,
\]
where we use Chebyshev's sum inequality to get 
\[
\sum_{j=2}^K\frac{1-\lambda_j}{\sqrt{\vartheta_j}}\leq \sum_{j=2}^K(1-\lambda_j)\frac{1}{K-1}\sum_{k=2}^K\sqrt{\frac{2NW}{k}}\leq 4\sum_{j=2}^K(1-\lambda_j)\,, 
\]
Lemma \ref{karnik_bound22} to control $\sum_{j=2}^K(1-\lambda_j)=O(N^{-2})$ since $c>2$.
For $I_3$, similarly we use Lemma \ref{taper_upper bound in frequency domain at 0, post-plunge} to control $U_j(0)$ and get
\begin{align*}
|I_3|\leq \frac{2}{W}\sum_{j=K+1}^{2NW+G}\Big(\sqrt{\frac{2NW}{j}}\vee 1\Big)\lambda_j^2 +\frac{2}{W}\sum_{j=2NW+G+1}^{N-1}\lambda_j^2 \,. 
\end{align*}
For the first summation, we use $\sqrt{\frac{2NW}{k+1}}=1+O\big(\frac{G}{NW}\big)$ and there are at most $2G$ nontrivial eigenvalues bounded by $1$, and hence the first summation is $O\big(\frac{\log(NW)\log(N)}{W}\big)$. For the second term, Lemma \ref{karnik_bound22} leads to $ \sum_{j=2NW+G+1}^{N-1}\lambda_k^2=O(N^{-1})$, which is dominated by the first term. 
As a result, we have
\begin{align*}
|I_3|=O\Big(\frac{\log(NW)\log(N)}{W}\Big)\,.
\end{align*}
\end{proof}

\begin{lemma}[When $f$ is close to $0$]\label{lemma: sum Uk(0)Uk(f) When f is close to 0}
Suppose $W=W(N)\to 0$ and $WN\to \infty$ as $N\to \infty$, and $K=\lfloor 2NW\rfloor-G$, where $G:=\lfloor \frac{2c}{\pi^2}\log(100NW+25)\log(N)\rfloor$ with $c>2$. 
When $N$ is sufficiently large and $\gamma>1$, 
\[
\Big|\sum_{j=0}^K U_{j}(0)U_{j}\Big(\frac{1}{\gamma N}\Big)\Big|=\frac{\gamma }{\pi}\sin\Big(\frac{\pi}{\gamma }\Big) N+O\Big(\frac{\log(NW)\log(N)}{W}\Big)  \,.
\]
Moreover, for $k=1,2,\ldots, \lfloor \frac{NW}{\log(NW)}\rfloor$, we have
\[
\Big|\sum_{j=0}^K U_{j}(0)U_{j}\Big(\frac{k}{N}\Big)\Big|=O\Big(\frac{\log(NW)\log(N)}{W}\Big)   \,,
\]
and 
\[
\Big|\sum_{j=0}^K U_{j}(0)U_{j}(f)\Big|\leq \frac{N}{\pi k}+O\Big(\frac{\log(NW)\log(N)}{W}\Big) 
\]
when $f\in \big(\frac{k}{N},\frac{k+1}{N}\big)$.
\end{lemma}

Note that $U_j(f)$ is the convolution of $U_j(k/N)$ and Dirichlet kernel; that is, $U_j(f)=\frac{1}{N}\sum_{k=0}^{N-1}U_j(k/N)D_N\big(f-\frac{k}{N}\big)$.
When $\gamma $ is large, by Taylor expansion, we have $\frac{\gamma}{\pi}\sin\big(\frac{\pi}{\gamma }\big)=1-\frac{\pi^2}{6\gamma ^2}+O(\gamma ^{-4})$, which quantifies the decay behavior of the main lobe centered at $0$.

\begin{proof}
The proof is similar to that for Lemma \ref{lemma summation of V_j(0)^2 is N}. Consider \eqref{EQ:proof of lemma of sumUkUk(f) decomposition when f <W}.
$I_0=D_N(f)=\frac{\sin(\pi N f)}{\sin(\pi f)}$.
By a direct expansion, we have
\[
|I_0|=\Big|D_N\Big(\frac{1}{\gamma N}\Big)\Big| = \frac{\gamma }{\pi}\sin\Big(\frac{\pi}{\gamma }\Big) N\Big(1+O\Big(\frac{1}{\gamma N}\Big)\Big)
\]
and when $f=k/N$ for $k=1,2,\ldots, \lfloor\sqrt{NW}/\pi\rfloor$, $I_0=0$. 
When $f\in \big(\frac{k}{N},\frac{k+1}{N}\big)$, 
\[
|I_0|=|D_N(f)|\leq \frac{\sin(\pi N f)}{2 f} \leq \frac{N}{\pi k} \,,
\]
Note that the maximization happens slightly before $k+1/2$.

For $I_1$, we apply a simple bound. By Cauchy-Schwartz, we have
\begin{align*}
\Big|\int_{|x|>W}D_N(f-x)D_N(x)dx\Big|\leq \sqrt{\int_{|x|>W}|D_N(f-x)|^2dx}\sqrt{\int_{|x|>W}|D_N(x)|^2dx}\,.
\end{align*}
For the first integration, note that by assumption $WN\to \infty$ and $f/ W\to 0$ uniformly when $N\to \infty$. Using $|D_N(t)|\leq \frac{1}{|\sin(\pi t)|} \leq \frac{1}{2|t|}$ since $\sin(x)\geq \frac{2x}{\pi}$ over $[0,\pi/2]$, we get 
\[
\int_{W}^{1/2}|D_N(f-x)|^2dx\leq \frac{1}{4(W-f)}-\frac{1}{2-4f}=\frac{1}{4W}\Big(1+O\Big(\frac{1}{\log(NW)}\Big)\Big)
\] 
and similarly for $\int^{-W}_{-1/2}|D_N(f-x)|^2dx=\frac{1}{4W}\big(1+O\Big(\frac{1}{\log(NW)}\Big)\big)$. 
The second integration is bounded directly by $2\int_W^{1/2}\frac{\pi^2}{4x^2}dx< \frac{\pi^2}{2W}$ using $\sin(x)\geq \frac{2x}{\pi}$ over $[0,\pi/2]$. As a result, 
\[
|I_1|=O(W^{-1})\,.
\]

For $I_2$, we use \eqref{control of Uk0 using pointwise control existing results} to control $U_0(0)$ and the fact that $0$ is the global maximum of $U_0$, and Lemma \ref{taper_upper bound in frequency domain at 0, post-plunge} to control $U_j(0)$ when $2\leq j\leq K$, and the rough bound \eqref{Uk worst case bound} to control $|U_j(f)|$. These lead to
\begin{align*}
|I_2|\,&\leq \frac{N^{1/4}}{W^{1/4}}(1-\lambda_0)|U_0(f)|+\frac{2}{\sqrt{W}} \sum_{j=2}^{K}(1-\lambda_k) \Big(\frac{2NW}{j}\Big)^{1/4}|U_j(f)|\\
&\leq \frac{N^{1/2}}{W^{1/2}}(1-\lambda_0)+C\frac{N^{3/8}}{W^{7/8}}\sum_{j=2}^K \frac{(1-\lambda_k)(j+1)^{1/8}}{j^{1/4}}\,,
\end{align*}
where $C=\frac{3^{9/4}}{2^{17/4}\pi^{1/8}}$. By Theorem \ref{karnik_bound}, we have $\frac{N^{1/2}}{W^{1/2}}(1-\lambda_0)=O(N^{-1})$.
Since $1-\lambda_j$ is increasing and $\frac{(j+1)^{1/8}}{j^{1/4}}$ is decreasing, we apply Chebyshev's sum inequality to bound
\[
\sum_{j=2}^K \frac{(1-\lambda_k)(j+1)^{1/8}}{j^{1/4}}\leq \sum_{j=2}^K(1-\lambda_j) \frac{1}{K-1}\sum_{k=2}^K \frac{(k+1)^{1/8}}{k^{1/4}} \leq 2N^{-2}K^{-1/8}\,,
\]
where the last bound comes from Lemma \ref{karnik_bound22} since $c>2$. Putting all together, we have
\begin{align*}
|I_2| =O(N^{-1})\,.
\end{align*}

For $I_3$, we have
\begin{align}
&|I_3|\leq \sum_{j=K+1}^{2NW+G}\lambda_{j}|U_{j}(0)||U_{j}(f)|+\sum_{j=2NW+G+1}^{N/2}\lambda_{j}|U_{j}(0)||U_{j}(f)|\,.\label{EQ:proof of lemma of sumUkUk(f) decomposition when f <W I3term}
\end{align}
Since $K$ is close to $2NW$, by Lemma \ref{taper_upper bound in frequency domain, control of 1-A/B-A}, when $k<2NW$, $W-w_k=O(N^{-1}\log(NW)\log(N))$, which means $w_k>\frac{W}{\log(NW)}$. Therefore, we can safely consider Theorem \ref{theorem:said2023non:summary}(d) to bound $U_k(f)$:
\[
|U_k(f)|^2\leq  |U_k(0)|^2\sqrt{\frac{1-A}{\cos(2\pi f)-A}\frac{1-B_k}{\cos(2\pi f)-B_k}}\,.
\] 
Since $W\to 0$ and  $f\leq \frac{W}{\log(NW)}$, which is much smaller than $W$, $\cos(2\pi f)-A=2\pi^2W^2\big(1-\frac{1}{\log(NW)^2}\big)(1+O(W^2))$ and $1-A=2\sin(\pi W)^2=2\pi^2W^2(1+O(W^2))$, and hence $\frac{1-A}{\cos(2\pi f)-A}=1+O(\frac{1}{\log(NW)^2})$.
Applying a trivial bound $\frac{1-B_k}{\cos(2\pi f)-B_k}<\sqrt{2}\frac{1-A}{\cos(2\pi f)-A}$ over $[0,\frac{W}{\log(NW)}]$ when $N$ is sufficiently large, we have
\[
|U_k(f)|\leq \sqrt{ 2}|U_k(0)|\,.
\]
The problem is then reduced to the same control of $I_3$ in the proof of Lemma \ref{lemma summation of V_j(0)^2 is N}, and we obtain
\[
|I_3|=O\Big(\frac{\log(NW)\log(N)}{W}\Big)\,.
\]

\end{proof}

\begin{lemma}[$f$ is away from $0$ but smaller than or equal to $W$]
Suppose $W=W(N)\to 0$ and $WN\to \infty$ as $N\to \infty$, and $K=\lfloor 2NW\rfloor-G$, where $G:=\lfloor \frac{2c}{\pi^2}\log(100NW+25)\log(N)\rfloor$ with $c>2$.  Denote $f_*:=\Big(1- \frac{G}{(NW)^{1/4}}\Big)W$. 
When $N$ is sufficiently large, we have
\[
\Big|\sum_{j=0}^K U_{j}(0)U_{j}(f)\Big|=O\Big(\frac{N^{1/8}\sqrt{\log(NW)\log(N)}}{W^{7/8}}\Big)\,,
\]
uniformly over $[\frac{W}{\log(NW)}, f_*]$, and  
\[
\Big|\sum_{j=0}^K U_{j}(0)U_{j}(f)\Big|=O\Big(\frac{N^{1/8}\log(NW)\log(N)}{W^{9/8}}\Big)
\]
uniformly over $(f_*,W]$.
\end{lemma}

\begin{proof} 
Again, consider \eqref{EQ:proof of lemma of sumUkUk(f) decomposition when f <W}.  $I_0=\int_{-1/2}^{1/2}D_N(f-x)D_N(x)dx=D_N(f)=\frac{\sin((2N+1)\pi f)}{\sin(\pi f)}$. Since $\big|\sum_{j=0}^K U_{j}(0)U_{j}\big(\frac{W}{\log(NW)}\big)\big|$ is already dominated by $\frac{\log(NW)\log(N)}{W}$, we simply bound $\sin((2N+1)\pi f)$ by $1$ and use $\sin(\pi f)\geq 2f$ to control 
\[
|I_0|\leq \frac{1}{2f}\leq \frac{\log(NW)}{2W}\,.
\]
For $I_1$, we apply the traditional $L^1$ bound of Dirichlet integral and have
\[
|I_1|\leq \frac{1}{\sin(\pi W)}\int_{-1/2}^{1/2}|D_N(f-x)|dx=O\Big(\frac{\log(N)}{W}\Big)\,.
\]
$I_2$ is controlled in the same way of that of Lemma \ref{lemma: sum Uk(0)Uk(f) When f is close to 0}, and we get $|I_2|=O(N^{-1})$.  $I_3$ is the main term we need to control. Decompose it into three terms:
\[
\underbrace{\sum_{j=2NW-G}^{2NW-1}\lambda_{j}U_{j}(0)U_{j}(f)}_{:=I_{3,1}}+\underbrace{\sum_{j=2NW}^{2NW+G}\lambda_{j}U_{j}(0)U_{j}(f)}_{:=I_{3,2}}+\underbrace{\sum_{j=2NW+G+1}^{N-1}\lambda_{j}U_{j}(0)U_{j}(f)}_{:=I_{3,3}}\,.
\] 
For $I_{3,1}$ and $I_{3,2}$, we control it on two regions. To determine the two regions, note that $w_k<W$ is increasing when $k$ increases from $1$ to $2NW$. With \eqref{proof of Ukf decay, wk and W in Bk and A} and the control in Lemma \ref{taper_upper bound in frequency domain, control of 1-A/B-A}, we know when $k\geq 2NW-G$, 
\[
w_k\geq w_{2NW-G}> \Big(1- \frac{0.591G}{4 NW}\Big)W>\Big(1- \frac{G}{(NW)^{1/4}}\Big)W:=f_*\,.
\]
Clearly, $f_*<W$ and 
\[
w_k-f_*>\frac{G}{(NW)^{1/4}}W\,. 
\]
Consider two regions, $[\frac{W}{\log(NW)},f_*]$ and $(f_*,W]$. 
First, we control over the region $[\frac{W}{\log(NW)},f_*]$. By Theorem \ref{theorem:said2023non:summary}(d), when $f\leq f_*<w_k$,  
\[
|U_k(f)|^2\leq  |U_k(0)|^2\sqrt{\frac{1-A}{\cos(2\pi f)-A}\frac{1-B_k}{\cos(2\pi f)-B_k}}\,.
\] 
We have $1-A=2\sin(\pi W)^2=2\pi^2W^2(1+O(W^2))$ and $\cos(2\pi f)-A=2\pi^2(W^2-f^2)(1+O(W^2))$, so on $[\frac{W}{\log(NW)},f_*]$,
\begin{align*}
\frac{1-A}{\cos(2\pi f)-A}\,&=\frac{W^2}{W^2-f^2}(1+O(W^2))\leq \frac{W^2}{W^2-f_*^2}(1+O(W^2)) \leq \frac{(NW)^{1/4}}{G}(1+O(W^2))\,.
\end{align*}
On the other hand,  on $[\frac{W}{\log(NW)},f_*]$, $\cos(2\pi f)-B_k=2\pi^2(w_k^2-f^2)(1+O(W^2))\geq 2\pi^2(w_k^2-f_*^2)(1+O(W^2)) \geq 4\pi^2w_k(w_k-f_*)(1+O(W^2)) \geq 4\pi^2 \frac{W^{3/4}}{N^{1/4}} Gw_k(1+O(W^2))$, 
and hence
\[
\frac{1-B_k}{\cos(2\pi f)-B_k}\leq\frac{(NW)^{1/4}}{2\pi G}(1+O(W^2))\,.
\] 
when $N$ is sufficiently large. Hence, 
\[
|U_k(f)|^2\leq |U_k(0)|^2\frac{(NW)^{1/4}}{2 G}
\] 
over $f\in[\frac{W}{\log(NW)},f_*]$, and with Lemma \ref{taper_upper bound in frequency domain at 0, post-plunge}, we have
\[
\lambda_k|U_k(0)||U_k(f)|\leq  \frac{(NW)^{1/8}}{\sqrt{2G}}|U_k(0)|^2 =O\Big(\frac{N^{1/8}}{W^{7/8}\sqrt{\log(NW)\log(N)}}\Big)\,,
\]
where we use $\lambda_{k}\leq 1$, and hence 
\[
|I_{3,1}|=O\Big(\frac{N^{1/8}\sqrt{\log(NW)\log(N)}}{W^{7/8}}\Big)\,.
\]
For $I_{3,2}$, we apply Lemma \ref{taper_upper bound in frequency domain, f<W k>2NW case} to get 
\[
|U_{j}(f)|\leq  \sqrt{\frac{\lambda_j}{W}}\Big(1+O\Big(\frac{\log(NW)\log(N)}{(NW)^{1/2}}\Big)\Big)
\] 
for $ 2NW\leq j\leq 2NW+G$ over $f\in[\frac{W}{\log(NW)},f_*]$ and Lemma \ref{taper_upper bound in frequency domain at 0, post-plunge} to get $U_{j}(0)\leq 2\sqrt{\frac{\lambda_{j}}{W}}$, and obtain 
\[
|I_{3,2}|\leq \frac{3}{W}\sum_{j=2NW+1}^{2NW+G} \lambda^2_j =O\Big( \frac{\log(NW)\log(N)}{W} \Big)
\]
when $f\in[\frac{W}{\log(NW)},f_*]$, where we use the fact that $\sum_{j=2NW+1}^{2NW+G}\lambda^2_{j}=O(\log(NW)\log(N))$. 

Over $[f_*,W]$, 
we apply the rough bound \eqref{Uk worst case bound} to control $|U_j(f)|$ and the fact that $\lambda_{2j}\leq 1$ to obtain
\[
|I_{3,1}+I_{3,2}|\leq \frac{3^{9/4}}{2^{11/4}\pi^{1/8}}\sum_{j=2NW-G}^{2NW+G} \sqrt{\frac{2\lambda_j}{W}}\frac{N^{1/8}}{W^{5/8}}=O\Big(\frac{N^{1/8}\log(NW)\log(N)}{W^{9/8}}\Big)\,.
\]

For $I_{3,3}$, we simply apply the rough bound \eqref{Uk worst case bound} to control the size of post-plunge $|U_j(f)|$, Lemma \ref{taper_upper bound in frequency domain at 0, post-plunge} to control $U_{j}(0)\leq2\sqrt{\frac{\lambda_{j}}{W}}$, and Chebychev's sum inequality to bound 
\[
\sum_{k=2NW+G+1}^{N-1}\!\!\!\!\!\!(k+1)^{1/8}\lambda_k^{3/2}\leq \frac{\sum_{k=2NW+G+1}^{N-1}(k+1)^{1/8}}{N-1-2NW-G}\!\!\!\!\!\!\sum_{j=2NW+G+1}^{N-1}\!\!\!\!\!\!\lambda_j^{3/2}\leq N^{1/8}\!\!\!\!\!\!\sum_{j=2NW+G+1}^{N-1}\!\!\!\!\!\!\lambda_j^{3/2}\,. 
\]
By Lemma \ref{karnik_bound22}, $\sum_{j=2NW+G+1}^{N-1}\lambda^{3/2}_{j}=O(N^{-3})$, since $c>2$. As a result, when $N$ is sufficiently large, we have  
\[
|I_{3,3}|\leq \frac{3^{7/4}}{2^{9/4}\pi^{1/8}}\frac{N^{1/4}}{W^{9/8}} \sum_{j=2NW+G+1}^{N-1}\lambda_j^{3/2}=O(N^{-1})\,.
\]

By putting all together, $|I_{3,1}|$ dominates $|I_0|+|I_1|+|I_2|+|I_{3,2}|+|I_{3,3}|$ over $f\in[\frac{W}{\log(NW)},f_*]$, so $\big|\sum_{j=0}^K U_{j}(0)U_{j}(f)\big|=O\big(\frac{N^{1/8}\sqrt{\log(NW)\log(N)}}{W^{7/8}}\big)$.
 Over $[f_*,W]$, $|I_{3,1}+I_{3,2}|$ dominates $|I_0|+|I_1|+|I_2|+|I_{3,3}|$, and we have $\big|\sum_{j=0}^K U_{j}(0)U_{j}(f)\big|=O\big(\frac{N^{1/8}\log(NW)\log(N)}{W^{9/8}}\big)$.
\end{proof}

\begin{lemma}[$f$ is in the outside band region]
Suppose $W=W(N)\to 0$ and $WN\to \infty$ as $N\to \infty$, and $K=\lfloor 2NW\rfloor-G$, where $G:=\lfloor \frac{4c}{\pi^2}\log(100NW+25)\log(N)\rfloor$ with $c>3$. Fix a small $\tau>0$. 
When $N$ is sufficiently large and $f\in (W,(1+\tau)W]$, we have
\[
\Big|\sum_{j=0}^K U_{2j}(0)U_{2j}(f)\Big|=O(N^{-1}W^{-1/2})\,;
\]
when  $f\in \big((1+\tau)W,\frac{1}{2}\big]$,
\[
\Big|\sum_{j=0}^K U_{2j}(0)U_{2j}(f)\Big|=O(N^{-3}W^{-1/2})\,.
\]
\end{lemma}

\begin{proof}
We use \eqref{control of Uk0 using pointwise control existing results} to control $U_0(0)$.
When $f\in (W,(1+\tau)W]$, by Lemmas \ref{taper_upper bound in frequency domain at 0, post-plunge} and the worst case control of $|U_{2j}(f)|$ by \eqref{bound Uk over W,1/2 for preplung2} in Lemma \ref{taper_upper bound in frequency domain, f>W k<2NW case}, we have
\begin{align*}
\sum_{j=0}^K |U_{2j}(0)U_{2j}(f)| \,&\leq 
\frac{N^{1/2}}{W^{1/2}}|U_0(f)|+\frac{2}{\sqrt{W}} \sum_{j=2}^K \Big(\frac{2NW}{j}\Big)^{1/4}  |U_{2j}(f)|\\
&\leq 2\frac{N^{5/2}}{W^{1/2}}\sqrt{1-\lambda_0}+ \frac{4 N^2}{\sqrt{W}} \sum_{j=2}^K \Big(\frac{2NW}{j}\Big)^{1/4} \sqrt{1 - \lambda_j} \,,
\end{align*}
Similarly, the first term is controlled by Theorem \ref{karnik_bound}. By Chebyshev's sum inequality, the second term is bounded by
\[
\sum_{j=2}^K \Big(\frac{2NW}{j}\Big)^{1/4} \sqrt{1 - \lambda_j}\leq \frac{1}{K-1}\sum_{k=2}^K \Big(\frac{2NW}{k}\Big)^{1/4}\sum_{j=2}^K \sqrt{1 - \lambda_j}\leq 2N^{-3}\,,
\]
where we apply Lemma \ref{karnik_bound22} in the last inequality since $c>3$.
Putting all together, we have 
\[
\Big|\sum_{j=0}^K U_{2j}(0)U_{2j}(f)\Big|=O(N^{-1}W^{-1/2})\,,
\]
When $f\in \big((1+\tau)W,\frac{1}{2}\big]$, by \eqref{bound Uk over W,1/2 for preplung} and the same proof, we get the claim.

\end{proof}

With the above lemmas, we immediately obtain the proof of Theorem \ref{theorem: bound Vf energy leakage uniform version} since $c>4$.

\bibliographystyle{plain}
\bibliography{refs}

\end{document}